\documentclass[12pt,a4paper]{article}

\usepackage[english]{babel}
\usepackage{epsfig}
\usepackage{amsmath}
\usepackage{amssymb}
\usepackage{amsbsy}
\usepackage{amsfonts}
\usepackage{mathrsfs} 
\usepackage{array}
\usepackage{verbatim}
\usepackage{graphicx}
\usepackage{color}
\usepackage{tabularx}
\usepackage{bm}
\usepackage{booktabs}
\usepackage[font=small]{caption}
\usepackage{subfig}
\usepackage{stmaryrd}
\usepackage{textcomp}

\usepackage{array}
\usepackage{float}

\usepackage{amsmath}
\usepackage{amssymb}
\usepackage{amsthm}

\usepackage{graphicx}
\usepackage{mdframed} 

\usepackage{multirow}
\usepackage{stackengine}
\usepackage{booktabs}
\usepackage{tabularx}

\usepackage{natbib}
\usepackage{url}

\usepackage{microtype}
\usepackage[title]{appendix}

\usepackage{tikz} 
\usetikzlibrary{cd}
\usetikzlibrary{positioning}
\usetikzlibrary{decorations.pathreplacing} 
\usetikzlibrary{patterns}

\newcommand{\klein}[1]{{\color{black}{#1}}}

\pgfdeclarepatternformonly{south west lines}{\pgfqpoint{-0pt}{-0pt}}{\pgfqpoint{3pt}{3pt}}{\pgfqpoint{3pt}{3pt}}{
    \pgfsetlinewidth{0.4pt}
    \pgfpathmoveto{\pgfqpoint{0pt}{0pt}}
    \pgfpathlineto{\pgfqpoint{3pt}{3pt}}
    \pgfpathmoveto{\pgfqpoint{2.8pt}{-.2pt}}
    \pgfpathlineto{\pgfqpoint{3.2pt}{.2pt}}
    \pgfpathmoveto{\pgfqpoint{-.2pt}{2.8pt}}
    \pgfpathlineto{\pgfqpoint{.2pt}{3.2pt}}
    \pgfusepath{stroke}}

\theoremstyle{definition}

\newcommand{\revision}[1]{{\color{black}{#1}}}

\begin{document}

\title{A one-step blended soundproof-compressible model with balanced data assimilation: theory and idealised tests}
\author{Ray Chew$^{(1)}$, Tommaso Benacchio$^{(2)}$,\\ Gottfried Hastermann$^{(1)}$, Rupert Klein$^{(1)}$}
\date{}
\maketitle

\begin{center}

{\small
$^{(1)}$ FB Mathematik \& Informatik \\
Freie Universit\"at Berlin \\
Arnimallee 6, 14195 Berlin, Germany\\
{\tt \{ray.chew,gottfried.hastermann\}@fu-berlin.de}
{\tt rupert.klein@math.fu-berlin.de}
}
{\small
\vskip 0.5cm
$^{(2)}$ MOX -- Modelling and Scientific Computing \\
Dipartimento di Matematica, Politecnico di Milano \\
Piazza Leonardo da Vinci 32, 20133 Milano, Italy\\
{\tt tommaso.benacchio@polimi.it}
}
\end{center}

\noindent
{\bf Keywords}:   semi-implicit models; finite volume methods; hyperbolic equations; compressible flow;
soundproof models; ensemble data assimilation methods.

\vspace*{0.5cm}

\noindent
{\bf AMS Subject Classification}:   65M08, 65Z99, 76M12,  76R99, 86A10

\vfill
{\noindent\footnotesize This work has been submitted to the Monthly Weather Review. Copyright in this work may be transferred without further notice.}

\pagebreak

\abstract{A challenge arising from the local Bayesian assimilation of data in an atmospheric flow simulation is the imbalances it may introduce. Acoustic fast-mode imbalances of the order of the slower dynamics can be negated by employing a blended numerical model with seamless access to the compressible and the soundproof pseudo-incompressible dynamics. Here, the blended modelling strategy by Benacchio et al., MWR, vol.~142~(2014) is upgraded in an advanced numerical framework and extended with a Bayesian local ensemble data assimilation method. Upon assimilation of data, the model configuration is switched to the pseudo-incompressible regime for one time-step. After that, the model configuration is switched back to the compressible model for the duration of the assimilation window. The switching between model regimes is repeated for each subsequent assimilation window. An improved blending strategy for the numerical model ensures that a single time-step in the pseudo-incompressible regime is sufficient to suppress imbalances coming from the initialisation and data assimilation. This improvement is based on three innovations: (i) the association of pressure fields computed at different stages of the numerical integration with actual time levels; (ii) a conversion of pressure-related variables between the model regimes derived from low Mach number asymptotics; and (iii) a judicious selection of the pressure variables used in converting numerical model states when a switch of models occurs. Idealised two-dimensional travelling vortex and buoyancy-driven bubble convection experiments show that acoustic imbalances arising from data assimilation can be eliminated by using this blended model, thereby achieving balanced analysis fields.}

\vfill

\paragraph{Significance statement}
Weather forecasting models use a combination of physics-based algorithms and meteorological measurements. A problem with combining outputs from the model with measurements of the atmosphere is that insignificant signals may generate noise and compromise the physical soundness of weather-relevant processes. By selecting atmospheric processes through the toggling of parameters in a mixed model, we propose to suppress the undesirable signals in an efficient way and retain the physical features of solutions produced by the model. The approach is validated here for acoustic imbalances using a compressible/pseudo-incompressible model pair. This development has the potential to improve the techniques used to bring observations into models and with them the quality of atmospheric model output.

\pagebreak

\section{Introduction}
\subsection{Motivation}
Dynamical processes in the atmosphere evolve on a range of spatio-temporal scales, most comprehensively expressed by the full compressible flow equations. Limit regimes, derived from the full compressible flow equations by scale analysis and asymptotics, describe reduced dynamics, examples being the soundproof anelastic and pseudo-incompressible models traditionally used at small- to mesoscale, and the hydrostatic primitive equations at large to planetary scales \citep{pedlosky2013geophysical, vallis2017atmospheric, klein2010scale}.

To access the dynamics of the full compressible flow equations and of their limit regimes, separate numerical schemes can be developed for each of the limiting models. From a computational perspective, \revision{however,} the \revision{discrepancies between numerical solutions} of different equation sets \revision{obtained} by \revision{essentially} the same numerical scheme can be \revision{substantially} smaller than the \revision{discrepancies} associated with the solution of \revision{one and the same} equation set by different numerical schemes \citep{smolarkiewicz2008conservative,klein2009asymptotics}.

\klein{\citet{bok2014}, \citet{klein2014using} and, separately, \citet{smolarkiewicz2014} developed discretisation schemes} for the compressible equations that allow access to the pseudo-incompressible model within a single numerical framework, showing equivalent results of both configurations in small- to mesoscale tests {involving acoustically balanced flows}. The blended analytical and numerical framework{ in \cite{bok2014,klein2014using}, within which the compressible to pseudo-incompressible transition is realised as a continuum of models controlled by an appropriate \textit{blending parameter}, was {conceptually} extended in \citet{klein2016} to include access to hydrostatic models. \citet{bk2019} then proposed a numerical implementation and achieved equivalence of hydrostatic and nonhydrostatic model solutions on large scales in the absence of vertically propagating acoustic modes.}

Balanced data assimilation provides a key motivation for blended numerical models. A problem with local data assimilation is the imbalance that \klein{it may induce} \citep{lorenc2003potential}. As the assimilation procedure does not take heed of \klein{specific characteristics of a flow, such as conservation of mass, momentum, and energy, or of particular smoothness properties, the initial balance of a flow state may be destroyed by the assimilation procedure}, see \cite{neef2006four} and more specifically \cite{greybush2011balance, bannister2015balance} on the effects of localisation on balanced analysis fields.

Physically, local data assimilation in a compressible framework \revision{can introduce} imbalances through fast acoustic modes with velocity amplitudes that may be of the same order of magnitude as the velocities found in the slowly evolving balanced dynamics of interest, with potentially destructive effects on overall solution quality \citep{hoheneggerschar2007}. Judicious use of a blended soundproof-compressible model can be employed to counteract this effect. Imbalances inherent in the initial pressure fields can be effectively \klein{reduced} by solving the initial time-steps \klein{of a simulation} in the pseudo-incompressible regime \klein{so that, upon the subsequent transition to the compressible regime over several further time steps,} the pressure field is balanced with respect to the initial velocities and potential temperature fields \revision{\citep{bok2014, klein2014using}}.
\revision{
More specifically, the blending method leverages a discrete orthogonal projection onto the space of pseudo-incompressible solutions. Therefore, the blending scheme used in an ensemble data assimilation framework yields the ensemble of balanced solutions closest to the analysis ensemble, measured in a norm weighted by the mass-weighted potential temperature.}

By extension of this insight, when mounting data assimilation on the numerics, \revision{a projection of the solution onto} the soundproof pseudo-incompressible model can \revision{suppress} the fast acoustic modes arising from the assimilation procedure. After suppression of the fast modes, the remaining time-steps until the next assimilation procedure are solved with the compressible model. \revision{As  this method makes use of the different dynamics modelled by the compressible and soundproof equation sets, it fundamentally deviates from existing methods to handle initialisation problems  such as} the post-analysis digital filter \citep[DFI, e.g.,][]{lynch1992initialization} \revision{and the incremental analysis update \citep[IAU, ][]{bloom1996data}}. \revision{These techniques act as low-pass filters, and} repeated application of the filter may have undesirable effects on long-term dynamics \revision{\citep{houtekamer2016review, polavarapu2004relationship}}.

Balance was also shown to improve with the choice of localisation space \citep{kepert2009covariance} and by allowing observations outside of a localisation radius to relax to a climatological mean \citep{flowerdew2015towards}. \cite{hastermann2021} \revision{compared} the effects of \revision{the blending approach} \revision{with those of} the post-analysis penalty method in achieving balanced analysis fields for highly oscillatory systems \revision{and found comparable improvements for both methods in the case of nonlinear balance relations.} See also \citet{zupanski2009theoretical} \revision{and} \citet{houtekamer2016review} for reviews of balanced atmospheric data assimilation.

\subsection{Contributions}
This paper proposes\revision{ a dynamics-driven method to achieve }balanced data assimilation using a blended numerical framework with the following advances:
\begin{itemize}
    \item One-step blending of the pseudo-incompressible and compressible models by instantaneous switching\klein{. This is} achieved by (a) accounting for the \klein{fact that Exner pressure fields computed at comparable stages within a time step correspond to different time levels in the compressible and soundproof model}; (b) judiciously converting the \klein{thermodynamic variables between the compressible and soundproof models motivated by low Mach number asymptotic arguments; and (c) carefully selecting, based on (a) and (b), the pressure variables used in converting numerical model states at the blending time interfaces. One-step blending} is a sizeable improvement over \klein{\cite{bok2014}}, who needed several intermediate time-steps for the blending procedure.
    
    \item Exploitation of the blended framework for balanced ensemble data assimilation. \revision{We employ an untuned data assimilation scheme that is known to introduce imbalances.} After each assimilation of data, a single time-step in the pseudo-incompressible model configuration is used to \revision{suppress} the fast acoustic imbalances. The model configuration is then switched back to the compressible model. In the reported \revision{idealised} experiments, balanced analysis fields are obtained \revision{by combining data assimilation and blending}\revision{, thus verifying the ability of the blended model to handle imbalances consistently with the underlying compressible and soundproof dynamics.}
\end{itemize}

The effects of data assimilation and blending on balanced solutions are investigated in \revision{the two-dimensional  numerical experiments of a travelling vortex and of a rising thermal in a vertical slice} \revision{\citep[see][]{ kadioglu2008,MendezCaroll1994,klein2009asymptotics}}. \revision{For} these tests, unbalanced \revision{and untuned} data assimilation is shown \revision{here} to destroy solution quality, while the use of blending effectively recovers the structure of the solution as evaluated by comparison with runs without data assimilation. \revision{Moreover, with the balanced data assimilation procedure, the solution quality of the observed quantities is maintained or improved independently of the size of the localisation region, which is an important tunable parameter of many sequential data assimilation procedures. The imbalances introduced by data assimilation in these idealised test cases can be quantified by scale analysis.}

The paper is structured as follows. Section \ref{sec:da_intro} contains a brief introduction to data assimilation and the Kalman filters \klein{considered here}. Section \ref{sec:blended_nm} reviews the blended numerical framework. Section \ref{sec:blending_scheme} proposes the new blending scheme and section \ref{sec:numerical_results} details the results of numerical experiments. The effectiveness of the one-step blended soundproof-compressible scheme is investigated for balanced data initialisation in section \ref{subsec:effectiveness_blending} and its application towards balanced data assimilation in \ref{subsec:da_results}. Section \ref{sec:diss_conc} contains discussion and conclusion. 

\section{Data assimilation: a quick primer}
\label{sec:da_intro}
Data assimilation is used in numerical weather prediction to improve forecast\revision{ing}. \revision{Existing approaches include} 4D-Var, which optimises model states over a finite time horizon in the past before launching a new prediction, and \revision{sequential assimilation procedures, which assimilate the available observations at specific points in time}. Here we focus on the latter, which are more susceptible to the problem of imbalances addressed in this paper \revision{due the local nature of these methods and especially when the localisation is severe \citep{cohn1998assessing, mitchell2002ensemble}}.

\klein{Modern weather forecasting techniques aim to represent the uncertainty of a forecast by generating an ensemble of likely candidates of model states. Such an ensemble can be understood as an approximate representation of a probability distribution over model states. The task of \revision{sequential} data assimilation is then as follows. Suppose we are given the probabilistic weight of each ensemble member at a previous instance in time, \textit{i.e.}, at the beginning of the current simulation window, together with the forward simulation states of all ensemble members at the current time, \textit{i.e.}, at the end of the simulation window. Then the \textit{prior probability distribution} $\text{pdf}_{\text{prior}}$ is represented by the model states at the new time level together with their probabilistic weights inherited from the beginning of the simulation window. Now we are to readjust the current states or the probabilistic weights of the ensemble members, at fixed time, such that the resulting \textit{posterior probability distribution} $\text{pdf}_{\text{post}}$ best reflects the observations that have arrived during the simulation window.}

\klein{The connection between $\text{pdf}_{\text{prior}}$ and $\text{pdf}_{\text{post}}$ can be established in a Bayesian framework. \revision{For this purpose we assume no model error and $\mathbf{x}_{\text{truth}}$ to be a perfectly resolved, \textit{true} model state. We denote the observation operator by \(\mathcal{H}\).  Then, observations that have arrived during the simulation window are subject to Gaussian distributed noise \(\epsilon\) and satisfy 
\begin{equation}
        \mathbf{y}_{\text{obs}} = \mathcal{H}(\mathbf{x}_{\text{truth}}) + \epsilon.
\end{equation}
}
Now Bayes' theorem gives
\begin{equation}
    \text{pdf}_{\text{post}}(\mathbf{x})
    = \text{pdf}(\mathbf{x} \vert \mathbf{y}_{\rm obs}) = \frac{\text{pdf}(\mathbf{y}_{\rm obs} \vert \mathbf{x})}{\text{pdf}(\mathbf{y}_{\rm obs})} \text{pdf}_{\text{prior}}(\mathbf{x}).
\end{equation}
Here $\text{pdf}(\mathbf{x} \vert \mathbf{y}_{\rm obs})$ is the conditional probability of state $\mathbf{x}$ given the observations $\mathbf{y}_{\rm obs}$ and $\text{pdf}(\mathbf{y}_{\rm obs} \vert \mathbf{x})$ is the probability of observation $\mathbf{y}_{\rm obs}$ given the state $\mathbf{x}$. The right-hand side of this equation is computable given the information before the data assimilation step, noting that the best available estimate of $\text{pdf}(\mathbf{y}_{\rm obs})$ is the expectation of $\text{pdf}(\mathbf{y}_{\rm obs} \vert \mathbf{x})$ with respect to $\mathbf{x}$ under the prior probability distribution.}
See \cite{wikle2007bayesian, reich2013ensemble} for more details on Bayesian data assimilation.

\subsection*{The Kalman filters}
Kalman filters are a family of popular Bayesian-based data assimilation methods \citep{kalman1960new} \klein{that assumes Gaussian shape for all probability densities so that they can be fully characterised by their \revision{means and} covariance matrices}. Identifying the prior with the term \textit{forecast} ($
f$), and the posterior with the term \textit{analysis} ($a$), the Kalman filter is
\begin{align}
    \mathbf{x}^a &= \mathbf{x}^f +
    \revision{
    {\pmb{\mathsf{B}}\mathcal{H}^T}
    {( \mathcal{H} \pmb{\mathsf{B}} \mathcal{H}^T + \pmb{\mathsf{R}} )}^{-1}
    }
    (\mathbf{y}_{\rm obs} - \mathcal{H}(\mathbf{x}^f))
    \nonumber \\
    &= \mathbf{x}^f + \pmb{\mathsf{K}} (\mathbf{y}_{\rm obs} - \mathcal{H} (\mathbf{x}^f)),
\end{align}
where $\pmb{\mathsf{B}} \in \mathbb{R}^{m \times m}$ and $\pmb{\mathsf{R}} \in \mathbb{R}^{l \times l}$ are the covariance matrices associated with the forecast and observations, respectively. $\pmb{\mathsf{K}}$ is the Kalman gain, which rewards the forecast if $\pmb{\mathsf{B}} \ll \pmb{\mathsf{R}}$ and penalises it if $\pmb{\mathsf{R}} \ll \pmb{\mathsf{B}}$.

\revision{A class of \revision{Monte Carlo-based} Kalman filters, the} ensemble Kalman filters, avoid the problem of high dimensionality by approximating the underlying probability density functions through \klein{the empirical distributions given by an ensemble of individual simulation states} \revision{\citep{reich2015probabilistic}}. As a consequence, ensemble-based methods are \klein{often} computationally more efficient \klein{than any scheme that aims to explicitly describe entire probability density functions}.

Specifically, for an ensemble of size $K$, the ensemble forecast is $\{ \mathbf{x}^{f}_1, \dots , \mathbf{x}^{f}_K \}$ and the ensemble's parametric information \klein{specifying its probability distribution is} updated by
\begin{subequations}
\begin{align}
    {\bar{\mathbf{x}}}^a &= {\bar{\mathbf{x}}}^f + \pmb{\mathsf{K}}^{\text{ens}} (\mathbf{y}_{\rm obs} - \mathcal{H}({\bar{\mathbf{x}}}^f)), 
    \label{subeqn:ensKalman_distance}\\
    \pmb{\mathsf{P}}_{K}^a &= \pmb{\mathsf{P}}_{K}^f - \pmb{\mathsf{K}}^{\text{ens}} \revision{\mathcal{H}} \pmb{\mathsf{P}}_{K}^f, 
    \label{subeqn:ensKalman_covar}\\
    \pmb{\mathsf{K}}^{\text{ens}} &= \revision{
    {\pmb{\mathsf{P}}_K^f \mathcal{H}^T}{( \mathcal{H} \pmb{\mathsf{P}}_K^f \mathcal{H}^T + \pmb{\mathsf{R}})}^{-1}},
    \label{subeqn:ensKalman_gain}
\end{align}
\end{subequations}
where $\bar{\mathbf{x}}^{a/f}$ is the ensemble mean and $\pmb{\mathsf{P}}^{a/f}_K \in \mathbb{R}^{K \times K}$ is the covariance associated with the ensemble. 

A drawback to the ensemble Kalman filter is that the covariance is determined by the spread of the ensemble and is therefore typically underestimated. However, ensemble inflation can be applied by multiplying the ensemble covariance by a constant factor larger than 1. This increases the covariance in the direction of the ensemble spread \revision{\citep{anderson2007adaptive, van2015nonlinear}}.

This paper uses the local ensemble transform Kalman filter (LETKF) data assimilation method \citep{hunt2007} \revision{based on the ensemble square-root filter (ESRF)}. The LETKF localises the observation covariance in such a way that observations farther away from the grid point under analysis have less influence, tapering off to zero influence for observations outside of a prescribed observation radius. The algorithm for the LETKF is \klein{provided} in Appendix A.

Localisation prevents spurious correlations of faraway observations while potentially reducing the complexity of the problem by making the observation covariance \klein{matrix closer to} diagonal \revision{\citep{hamill2001distance, houtekamer1998data}. After localisation, the analysis is only performed on a smaller local region, and the global analysis ensemble comprises different linear combinations of the ensemble members in each of these local regions. This allows the ensemble to represent a higher-dimensional space than one constrained by the ensemble size \citep{fukumori2002partitioned, mitchell2002ensemble}. A smaller ensemble size may necessitate more severe localisation. 

\revision{
When applying the LETKF, there are two potential sources for imbalances. In the case of a nonlinear balance relation, the LETKF fails to recover the desired balance due to its local linear construction. Even without localisation and for a given observation, the analysis ensemble of the ESRF is obtained as a linear combination of the forecast ensemble. In the case of linear balances, the situation is more subtle. On one hand the ESRF is capable of resolving linear balances due to its linear construction. On the other hand the LETKF, utilising localisation, does not act as a linear map on the global fields and therefore does not necessarily preserve the balance relation.} \revision{Numerical experiments in this paper investigate imbalances arising from both these sources.}
}
A smooth localisation function, such as the truncated Gaussian function or the \cite{gaspari1999construction} function, may be used \revision{to keep the resulting field sufficiently smooth}.


\section{The blended numerical model}
\label{sec:blended_nm}

\subsection{Governing equations}
In a rotating three-dimensional Cartesian domain, the adiabatic, dry compressible fluid flow equations for an ideal gas under gravity are:
\begin{subequations}
\begin{align}
\rho_t + \nabla_{\Vert} \cdot \left( \rho \mathbf{u} \right) + (\rho w)_z & = 0,\label{eqn:governing_a}\\
\left( \rho \mathbf{u} \right)_t + \nabla_{\Vert} \cdot (\rho \mathbf{u} \circ \mathbf{u}) + (\rho w \mathbf{u})_z &= -\left[ c_p P \nabla_{\Vert} \pi + f\mathbf{k} \times \rho \mathbf{u} \right],\label{eqn:governing_b}\\
\left( \rho w \right)_t + \nabla_\Vert \cdot (\rho \mathbf{u} w) + (\rho w^2)_z  &= -\left( c_P P \pi_z + \rho g \right),\label{eqn:governing_d}\\
\alpha_P P_t + \nabla_{\Vert} \cdot \left( P \mathbf{u} \right) + (P w)_z &= 0,\label{eqn:governing_c}
\end{align}
\label{eqn:governing}%
\end{subequations}
where $\rho$ is the density, $\mathbf{u} = (u,v)$ the vector of horizontal velocities and $w$ the vertical velocity, $P$ is the mass-weighted potential temperature and $\pi$ is the Exner pressure. $f$ is the Coriolis parameter on the horizontal $(x,y)$-plane, $\mathbf{k}$ a unit vector in the vertical direction and $\times$ the cross product. $g$ is the acceleration of gravity acting in the direction of $\mathbf{k}$. $\circ$ denotes the tensor product, $\nabla_{\Vert}$ denotes the horizontal gradient while the subscripts $t$ and $z$ denote the partial derivatives with respect to time $t$ and the vertical coordinate $z$. \revision{$\pi$ and $P$} are related to the thermodynamic pressure $p$ by the equation of state,
\begin{align}
\pi = \left( \frac{p}{p_{\text{ref}}} \right)^{R/c_p},\qquad
P = \frac{p_{\text{ref}}}{R} \left( \frac{p}{p_\text{{ref}}} \right)^{c_v/c_p} = \rho \Theta,
\label{eqn:eos}
\end{align}
where $p_{\text{ref}}$ is a reference pressure, $c_p$ and $c_v$ are the heat capacities at constant pressure and constant volume, $R = c_p - c_v$ is the ideal gas constant, and $\Theta$ is the potential temperature. The parameter $\alpha_P$ tunes between the compressible and the pseudo-incompressible model \revision{\citep{durran1989,KleinEtAl2010}}.

Identifying $\chi$ with the inverse potential temperature
\begin{equation}
    \chi = \frac{1}{\Theta},
    \label{eqn:chi}
\end{equation}
the Exner pressure and inverse potential temperature can be decomposed as
\begin{subequations}
\begin{align}
    \pi &= \bar{\pi} + \pi^\prime \qquad \text{and}
    \label{eqn:pi_expansion} \\
    \chi &= \bar{\chi} + \chi^\prime,
    \label{eqn:chi_expansion}
\end{align}
\label{eqn:pi_chi_expansions}%
\end{subequations}
where the bar denotes a \klein{hydrostatic} background state, which depends only on the vertical coordinate, and the prime \revision{denotes} the perturbation.
\noindent Rewriting \eqref{eqn:governing} with \eqref{eqn:chi} yields
\begin{subequations}
\begin{align}
    \rho_t + \nabla_{\Vert} \cdot (P \mathbf{u} \chi) + (P w \chi)_z &= 0,
    \label{eqn:chi_governing_a} \\
    (\rho \mathbf{u})_t + \nabla_{\Vert} \cdot (P \mathbf{u} \circ \chi \mathbf{u}) + (P w \chi \mathbf{u})_z &= -\left[ c_p P \nabla_\Vert \pi + f \mathbf{k} \times \rho \mathbf{u} \right],
    \label{eqn:chi_governing_b} \\
   (\rho w)_t + \nabla_{\Vert} \cdot (P \mathbf{u} \chi w) + (P w \chi w)_z &= -\left( c_p P \pi_z + \rho g \right),
    \label{eqn:chi_governing_d} \\
    \alpha_P P_t + \nabla_\Vert \cdot (P \mathbf{u}) + (P w)_z &= 0.
    \label{eqn:chi_governing_c}
\end{align}
\label{eqn:chi_governing}%
\end{subequations}
Using the notation of \cite{smolarkiewicz2014} and \cite{bk2019},
\begin{equation}
    \Psi = (\chi , \chi \mathbf{u}, \chi w, \chi^\prime),
    \label{eqn:Psi_set}
\end{equation}
\eqref{eqn:chi_governing} can be written compactly as
\begin{subequations}
\begin{align}
    (P \Psi)_t + \mathcal{A}(\Psi ; P \mathbf{v}) &= Q(\Psi ; P),
    \label{eqn:compact_form_a}\\
    \alpha_P P_t + \nabla \cdot (P \mathbf{v}) &= 0,
    \label{eqn:compact_form_b}
\end{align}
\label{eqn:compact_form}%
\end{subequations}
where $\mathbf{v}=(u,v,w)$ subsumes the three-dimensional velocity fields, $\mathcal{A}(\Psi ; P \mathbf{v})$ denotes the advection of the quantity $\Psi$ given the advective fluxes $P \mathbf{v}$, while $Q(\Psi; P)$ describes the effect on the right-hand side of \eqref{eqn:chi_governing} on $\Psi$ given $P$.

From \eqref{eqn:eos}, $P$ is a function of $\pi$ \revision{only},
\begin{equation}
    P(\pi) = \frac{p_{\text{ref}}}{R} \pi^{\frac{1}{\gamma-1}},
    \label{eqn:P-pi_relation}
\end{equation}
where $\gamma = c_p / c_v$ is the isentropic exponent. With \eqref{eqn:P-pi_relation}, \eqref{eqn:compact_form_b} becomes
\begin{equation}
    \alpha_P \left( \frac{\partial P}{\partial \pi} \right) \pi_t = - \nabla \cdot (P \mathbf{v}).
    \label{eqn:dPdpi}
\end{equation}


\subsection{Summary of the numerical scheme}

Equation \eqref{eqn:chi_governing_c} is discretised in time with an implicit midpoint method,
\begin{equation}
    \alpha_P P^{n+1} = \alpha_P P^n - \Delta t ~ \nabla \cdot (P\mathbf{v})^{n+1/2}.
\end{equation}
In order to obtain the advective fluxes at the half time-level, the time-update for equations \eqref{eqn:compact_form} is split into advective and non-advective terms. The advection terms on the left are updated by
\begin{subequations}
\begin{align}
    (P \Psi)^\# &= \mathcal{A}^{\Delta t/2}_{\text{1st}} \left[ \Psi^n ; (P \mathbf{v}^n )\right],
    \label{eqn:half_time_lhs_a} \\
    \alpha_P P^\# &= \alpha_P P^n - \frac{\Delta t}{2} \tilde{\nabla} \cdot (P \mathbf{v})^n,
    \label{eqn:half_time_lhs_b}
\end{align}
\label{eqn:half_time_lhs}%
\end{subequations}
where $\tilde{\nabla}$ is the discrete divergence and $\mathcal{A}_{\text{1st}}$ is an advection scheme corresponding to the half time-level update. The terms on the right are then advanced using an implicit Euler method,
\begin{subequations}
\begin{align}
    (P \Psi)^{n+1/2} &= (P \Psi)^{\#} + \frac{\Delta t}{2} Q(\Psi^{n+1/2} ; P^{n+1/2}),
    \label{eqn:half_time_rhs_a} \\
    \alpha_P P^{n+1/2} &= \alpha_P P^n - \frac{\Delta t}{2} \tilde{\nabla} \cdot (P \mathbf{v})^{n+1/2}.
    \label{eqn:half_time_rhs_b}
\end{align}
\label{eqn:half_time_rhs}%
\end{subequations}
Expressions \eqref{eqn:half_time_lhs} and \eqref{eqn:half_time_rhs} yield the advective fluxes at the half time-level.

Subsequently, the quantities $\Psi$ are updated to the full time-level with an explicit Euler half step followed by a full advection step and a final implicit Euler half step,
\begin{subequations}
\begin{align}
    (P \Psi)^{*} &= (P \Psi)^{n} + \frac{\Delta t}{2} Q(\Psi^{n} ; P^{n}),
    \label{eqn:full_time_a} \\
    (P \Psi)^{**} &= \mathcal{A}^{\Delta t}_{\text{2nd}} \left[ \Psi^* ; (P \mathbf{v}^{n+1/2} )\right],
    \label{eqn:full_time_b} \\
    (P \Psi)^{n+1} &= (P \Psi)^{**} + \frac{\Delta t}{2} Q(\Psi^{n+1} ; P^{n+1}),
    \label{eqn:full_time_c} \\
    \alpha_P P^{n+1} &= \alpha_P P^n - \Delta t \tilde{\nabla} \cdot (P \mathbf{v})^{n+1/2},
    \label{eqn:full_time_d}
\end{align}
\label{eqn:full_time}%
\end{subequations}
yielding a second-order accurate one-step method \citep{bk2019, smolarkiewicz1991forward,smolarkiewicz1993forward}.

A first-order Runge-Kutta method is used for the advection operator $\mathcal{A}^{\Delta t/2}_{\text{1st}}$ in \eqref{eqn:half_time_lhs_a} while second-order Strang splitting is used for $\mathcal{A}^{\Delta t}_{\text{2nd}}$ in \eqref{eqn:full_time_b}. The former is necessary for the time-level analysis in section \ref{sec:blending_scheme} to hold and the latter \revision{guarantees second order in time of the overall scheme}. The spatial discretisation of the numerical scheme is based on a finite volume framework, for more details see section 4 in \cite{bk2019}.

\subsection{Pseudo-incompressible regime}
\klein{The switch $\alpha_P$ in \eqref{eqn:governing} toggles access to the pseudo-incompressible model \citep[$\alpha_P = 0$,][]{durran1989},  
\begin{subequations}
\begin{align}
\rho_t + \nabla_\Vert \cdot \left( \rho \mathbf{v} \right) + (\rho w)_z & = 0,\label{eqn:pi_a}\\
\left( \rho \mathbf{u} \right)_t + \nabla_{\Vert} \cdot (\rho \mathbf{u} \circ \mathbf{u}) + (\rho w \mathbf{u})_z &= -\left[ c_p P \nabla \pi + f\mathbf{k} \times \rho \mathbf{u} \right],\label{eqn:pi_b}\\
\left( \rho w \right)_t + \nabla_{\Vert} \cdot (\rho \mathbf{u} w) + (\rho w^2)_z &= -\left( c_p P \pi_z + \rho g \right),\label{eqn:pi_d}\\
\nabla_{\Vert} \cdot \left( P \mathbf{u} \right) + (P w)_z &= 0.\label{eqn:pi_c}
\end{align}
\label{eqn:pi}%
\end{subequations}
where \eqref{eqn:pi_c} enforces the soundproof divergence constraint. See  \cite{klein2009asymptotics,klein2012thermodynamic,klein2016} for details of this formulation.}

\section{Single time-step soundproof-\\compressible transition}
\label{sec:blending_scheme}
In the following, a conversion of pressure-related quantities, motivated by low Mach number asymptotics and applied prior to the model transitions, is proposed which allows for model switching within a single time-step.

\subsection{Time-level of the pressure-related variables}
\label{subsec:pi_time_level}
\klein{In the simpler non-rotating case without gravity ($g$, $f = 0$), the update for the momentum equation multiplied by the potential temperature $\Theta$ in \eqref{eqn:half_time_rhs_a} and \eqref{eqn:full_time_c} read}
\begin{equation}
(P \mathbf{v})_t = - c_p (P \Theta)^{\text{adv}} \nabla \pi,    
\label{eqn:general_rhs_update}
\end{equation}
where the superscript \textit{adv} denotes the quantity that becomes available after the advection substeps \eqref{eqn:half_time_lhs_a} and \eqref{eqn:full_time_b}. Applying an implicit Euler discretisation to \eqref{eqn:general_rhs_update}, \klein{we find}
\begin{equation}
    (P \mathbf{v})^{\text{out}} = (P \mathbf{v})^{\text{in}} - \delta t \, c_p (P \Theta)^{\text{adv}} \tilde{\nabla} \pi^{\text{out}},
    \label{eqn:discretised_rhs_update}
\end{equation}
where the superscript \textit{in} denotes the quantities at the time-level corresponding to the start of the time-step and \textit{out} at the end. $\delta t \leq \Delta t$ is an arbitrary time-step.

\subsubsection{The compressible equations}
\noindent For the case $\alpha_P = 1$, using \eqref{eqn:P-pi_relation} a discretisation of the left-hand side of \eqref{eqn:dPdpi} yields at the half time-level
\begin{equation}
    P^{n+1/2} - P^n = \left( \frac{\partial P}{\partial \pi} \right)^{\#} (\pi^{n+1/2} - \pi^n),
    \label{eqn:htl_dPdpi}
\end{equation}
where $\left( {\partial P} / {\partial \pi} \right)^{\#}$ is obtained from $P$ after the advection step at the half-time level, \eqref{eqn:half_time_lhs}. Substituting \eqref{eqn:htl_dPdpi} into \eqref{eqn:half_time_rhs_b}, 
\begin{equation}
    \left( \frac{\partial P}{\partial \pi} \right)^{\#} (\pi^{n+1/2} - \pi^n) = - \frac{\Delta t}{2} \tilde{\nabla} \cdot (P \mathbf{v})^{n+1/2}.
    \label{eqn:dPdpi_rhs_hts_update}
\end{equation}
\revision{$(P\mathbf{v})^{\rm in}$ in \eqref{eqn:discretised_rhs_update} is the solution of the advection terms in \eqref{eqn:half_time_lhs},
\begin{equation}
    (P\mathbf{v})^{\rm in} = (P \mathbf{v})^{n} - \frac{\Delta t}{2} \tilde{\nabla} \cdot (P \mathbf{v} \circ \mathbf{v})^n,
\end{equation}
with $\delta t = \Delta t / 2$. }Identifying \textit{out} with $n+1/2$ 
and rearranging, \eqref{eqn:dPdpi_rhs_hts_update} becomes
\begin{align}
    &\left( \frac{\partial P}{\partial \pi} \right)^{\#} \pi^{n+1/2} - \left(\frac{\Delta t}{2}\right)^2 \tilde{\nabla} \cdot \left[ \, c_p (P \Theta)^{\#} \tilde{\nabla} \pi^{n+1/2} \right] 
    = \nonumber \\
    &\qquad \left( \frac{\partial P}{\partial \pi} \right)^{\#} \pi^n - \frac{\Delta t}{2} \tilde{\nabla} \cdot (P\mathbf{v})^n + \left(\frac{\Delta t}{2}\right)^2 \tilde{\nabla} \cdot \tilde{\nabla} \cdot (P\mathbf{v} \circ \mathbf{v})^{n}.
    \label{eqn:comp_htl_fix}
\end{align}
which fixes the time-level of $\pi$ after the half-time step of \eqref{eqn:half_time_lhs} and \eqref{eqn:half_time_rhs} at $n+1/2$.

For the full-time stepping of \eqref{eqn:full_time}, a similar procedure yields
\revision{\begin{align}
    &\left( \frac{\partial P}{\partial \pi} \right)^{\#} \pi^{n+1} - \left(\frac{\Delta t}{2}\right)^2 \tilde{\nabla} \cdot \left[ \, c_p (P \Theta)^{\#} \tilde{\nabla} \pi^{n+1} \right] = 
    \nonumber \\
    &\qquad \quad \left( \frac{\partial P}{\partial \pi} \right)^{\#} \pi^n - \Delta t \, \tilde{\nabla} \cdot (P\mathbf{v})^n + \frac{(\Delta t)^2}{2} \tilde{\nabla} \cdot \tilde{\nabla} \cdot (P\mathbf{v} \circ \mathbf{v})^{n+1/2} 
    \nonumber \\
    &\qquad \qquad + \left(\frac{\Delta t}{2}\right)^2 \tilde{\nabla} \cdot \left[ \, c_p (P \Theta)^{\#} \tilde{\nabla} \pi^{n} \right].
    \label{eqn:comp_ftl_fix}
\end{align}}
From \eqref{eqn:comp_htl_fix}, $\pi$ is at time-level $n+1/2$ after the half-time stepping \eqref{eqn:half_time_rhs} while \eqref{eqn:comp_ftl_fix} starts with $\pi$ at time-level $n$ for the full-time stepping \eqref{eqn:full_time}. Therefore, the time-level of $\pi$ has to be reset from $n+1/2$ to $n$ after the half-time step \eqref{eqn:half_time_rhs} and before the full time-step \eqref{eqn:full_time}. Furthermore, the time-level of $\pi$ after the full time-step \eqref{eqn:full_time} is $n+1$ as intended.

\subsubsection{The pseudo-incompressible equations}
For $\alpha_P = 0$, the coupling between $P$ and $\pi$ in $\eqref{eqn:dPdpi}$ no longer holds and the two variables decouple, leading to
\begin{equation}
    \nabla \cdot (P \mathbf{v}) = 0.
    \label{eqn:divergence_constraint}
\end{equation}
\klein{Enforcing this divergence constraint for} the left-hand side of \eqref{eqn:discretised_rhs_update}, we obtain
\begin{equation}
    \tilde{\nabla} \cdot (P \mathbf{v})^{\text{in}} = \tilde{\nabla} \cdot \left( \delta t \, c_p (P \Theta)^{\text{adv}} \tilde{\nabla} \pi^{\text{out}} \right).
    \label{eqn:psinc_discretised_rhs_update}
\end{equation}
At the half-time level, $(P \mathbf{v})^{\text{in}}$ is the solution of \eqref{eqn:half_time_lhs} comprising the half time-step advection. Therefore,
\begin{equation}
    \tilde{\nabla} \cdot (P \mathbf{v})^{\text{in}} = \tilde{\nabla} \cdot \left[ (P \mathbf{v})^n + \frac{\Delta t}{2} \tilde{\partial_t} (P \mathbf{v})^{\#} \right],
    \label{eqn:intermediate_tla_pi}
\end{equation}
where $\tilde{\partial_t}$ is the discrete partial time derivative. Assuming that the divergence constraint \eqref{eqn:divergence_constraint} has been satisfied at the end of time-step $(n-1)$, the first term on the right-hand side vanishes. As the second term is \klein{generated by} \eqref{eqn:half_time_lhs} starting at time-level $n$, \textit{i.e.}, \klein{by an explicit} advection step \klein{associated with} the left-hand side of \eqref{eqn:chi_governing_b} and \eqref{eqn:chi_governing_d} multiplied by $\Theta$,
\begin{equation}
    \tilde{\partial_t} (P \mathbf{v})^\# + \tilde{\nabla} \cdot (P \mathbf{v} \cdot \mathbf{v})^n = 0,
\end{equation}
equation \eqref{eqn:intermediate_tla_pi} becomes
\begin{equation}
    \tilde{\nabla} \cdot (P \mathbf{v})^{\text{in}} = - \frac{\Delta t}{2} \tilde{\nabla} \cdot \left[ \tilde{\nabla} \cdot (P \mathbf{v} \circ \mathbf{v})^n\right].
    \label{eqn:psinc_pu_in}
\end{equation}
Inserting \eqref{eqn:psinc_pu_in} back into \eqref{eqn:psinc_discretised_rhs_update}, with $\delta t = \Delta t / 2$ and \textit{adv} as \textit{\#},
\begin{equation}
    \tilde{\nabla} \cdot \left( \frac{\Delta t}{2} c_p (P \Theta)^{\#} \tilde{\nabla} \pi^{n} \right) = - \frac{\Delta t}{2} \tilde{\nabla} \cdot \tilde{\nabla} \cdot (P \mathbf{v} \circ \mathbf{v})^n,
    \label{eqn:psinc_hts_fix}
\end{equation}
where the right-hand side has fixed the time-level of $\pi^{\text{out}}$ at $n$.

For the full time-stepping, $(P \mathbf{v})^{\text{in}}$ is the solution of \eqref{eqn:full_time_b} and so \eqref{eqn:psinc_discretised_rhs_update} is
\begin{equation}
    \tilde{\nabla} \cdot (P \mathbf{v})^{**} = \tilde{\nabla} \cdot \left( \frac{\Delta t}{2} c_p (P \Theta)^{**} \tilde{\nabla} \pi^{\text{out}} \right),
    \label{eqn:psinc_ftl_rhs_update}
\end{equation}
with
\begin{equation}
    \tilde{\nabla} \cdot (P \mathbf{v})^{**} = \tilde{\nabla} \cdot \left[ (P \mathbf{v})^n - \frac{\Delta t}{2} c_p (P \Theta)^{\#} \tilde{\nabla} \pi^n + \Delta t \, \tilde{\partial_t} (P \mathbf{v})^{**} \right],
    \label{eqn:psinc_full_expansion}
\end{equation}
where the second term in the square brackets is a correction of $(P \mathbf{v})$ from the implicit substep at the half time-level \eqref{eqn:half_time_rhs} and the third term is the solution of the advection substep at the full time-level. Assuming again that the divergence constraint \eqref{eqn:divergence_constraint} has been fulfilled at the onset, \begin{equation}
\tilde{\nabla} \cdot (P \mathbf{v})^n = 0.
\end{equation} Substitute \eqref{eqn:psinc_hts_fix} into \eqref{eqn:psinc_full_expansion},
\begin{equation}
    \tilde{\nabla} \cdot (P \mathbf{v})^{**} = \tilde{\nabla} \cdot \left[ \frac{\Delta t}{2} \tilde{\nabla} \cdot (P \mathbf{v} \circ \mathbf{v})^n + \Delta t \tilde{\partial_t} (P \mathbf{v})^{**} \right],
    \label{eqn:49_into_51}
\end{equation}
and note that advection substep \eqref{eqn:full_time_b} solves the left-hand side of \eqref{eqn:chi_governing_b} and \eqref{eqn:chi_governing_d} multiplied with $\Theta$, i.e,
\begin{equation}
    \tilde{\partial_t} (P \mathbf{v})^{**} + \tilde{\nabla} \cdot (P \mathbf{v} \circ \mathbf{v})^{n+1/2} = 0,
    \label{adv_left}
\end{equation}
where the half time-level of the second term emerges from the solution of substeps \eqref{eqn:full_time_a} and \eqref{eqn:full_time_b} under the advecting fluxes $(P\mathbf{v})^{n+1/2}$. Putting \eqref{eqn:49_into_51} and \eqref{adv_left} together yields, after re-elaborations,
\begin{equation}
    \tilde{\nabla} \cdot (P \mathbf{v})^{**} = - \frac{\Delta t}{2} \left[ 1 + \frac{\Delta t}{2} \tilde{\partial_t} \right] \tilde{\nabla} \cdot \left[ \tilde{\nabla} \cdot (P \mathbf{v} \circ \mathbf{v})^{n+1/2}\right ].
    \label{eqn:long_working}
\end{equation}
Inserting \eqref{eqn:long_working} back into \eqref{eqn:psinc_ftl_rhs_update} gives
\begin{equation}
    \tilde{\nabla} \cdot \left( c_p (P \Theta)^{**} \tilde{\nabla} \pi^{n+1} \right) = - \left( 1 + \frac{\Delta t}{2} \tilde{\partial_t} \right) \tilde{\nabla} \cdot \left[ \tilde{\nabla} \cdot (P \mathbf{v} \circ \mathbf{v})^{n+1/2}\right],
    \label{eqn:psinc_ftl_fix}
\end{equation}
fixing the time-level of $\pi^{\text{out}}$ at $n+1$, since the right-hand side is a half time-step advancement from the $n+1/2$ time-level.

\begin{figure}
\centering
\includegraphics[scale=1.0]{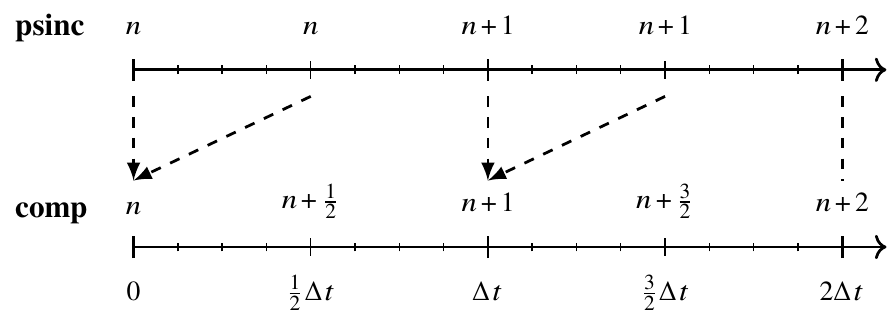}
\caption{Summary of the time-levels of $\pi$ for the pseudo-incompressible (\textbf{psinc}) and the compressible (\textbf{comp}) solutions in the numerical scheme. $\Delta t$ indicates time-step size. The dashed lines relate the $\pi$'s at the same time-level between the two models. The two valid choices of $\pi$ in the pseudo-incompressible to compressible blending, at time $t^n$ and $t^{n+1}$, are depicted with arrows.}
\label{fig:timeline}
\end{figure}

\klein{In contrast to} the compressible case, expressions \eqref{eqn:psinc_hts_fix} and \eqref{eqn:psinc_ftl_fix} imply that Exner pressure $\pi$ after the half-step $\eqref{eqn:half_time_lhs}$ and $\eqref{eqn:half_time_rhs}$ is at the time-level $n$, and \revision{could} be used as the input to $\eqref{eqn:full_time}$ \revision{as an alternative to using the Exner pressure obtained at the end of time step $n-1$. With this option, $\pi$ does not have to be} reset to time-level $n$ \revision{after the half-time predictor} for the pseudo-incompressible solve. Figure \ref{fig:timeline} summarises the time-level analysis of $\pi$.

\subsection{Conversion of the pressure-related variables}
\label{ssec:Conversion}

Expression \eqref{eqn:pi_expansion} separates the background Exner pressure from its perturbation. For \klein{low Mach number flows, $\mathrm{Ma} \ll 1$, such a separation is naturally induced by} the asymptotic expansion
\begin{equation}
\label{eqn:pi_expnsn}
    \pi = \pi^{(0)} + \mathrm{Ma}^2 \pi^{(1)} + \dots ~,
\end{equation}
where $\mathrm{Ma} = u_{\text{ref}} / c_{\text{ref}}$ for reference velocity $u_{\text{ref}}$ and speed of sound $c_{\text{ref}}$. Substituting this expansion into \eqref{eqn:dPdpi} \klein{yields}
\begin{equation}
\label{eqn:P_pi}
    \alpha_P \left( \frac{\partial P}{\partial \pi} \right) \mathrm{Ma}^2 \pi^{(1)}_t = - \nabla \cdot \left( P \mathbf{v} \right),
\end{equation}
based on which we can blend the pseudo-incompressible  ($\alpha_P=0$) and compressible ($\alpha_P=1$) models. Using \eqref{eqn:P-pi_relation}, the compressible $P_{\text{comp}}$ is obtained from the pseudo-incompressible model variables as
\begin{equation}
    P_{\text{comp}} = \left( P_{\text{psinc}}^{\gamma-1} + \mathrm{Ma}^2 \pi_{\text{psinc}} \right)^{\frac{1}{\gamma-1}}.
    \label{eqn:P_conv_psinc_to_comp}
\end{equation} 

By inverting \eqref{eqn:P_conv_psinc_to_comp}, the pseudo-incompressible $P_{\text{psinc}}$ is obtained as
\begin{equation}
  P_{\text{psinc}} = \left( P_{\text{comp}}^{\gamma-1} - \mathrm{Ma}^2 \pi_{\text{comp}} \right)^{\frac{1}{\gamma-1}}.
  \label{eqn:P_conv_comp_to_psinc}
\end{equation}
Therefore, at the blending time interfaces between the compressible and the pseudo-incompressible configurations, expression \eqref{eqn:P_conv_psinc_to_comp} or \eqref{eqn:P_conv_comp_to_psinc} is applied depending on the direction of the transition.

\subsection{Association of perturbation variables between the compressible and soundproof models}
\label{subsec:pi_choice}
The time-level analysis of $\pi$ in section \ref{subsec:pi_time_level} demonstrated that, in a pseudo-incompressible solve, both the Exner pressure solution after the full time-step from $t^n$ to $t^{n+1}$ and that obtained after the \textit{subsequent} half time-step are associated with the same time-level $t^{n+1}$. 

Consider then the compressible to pseudo-incompressible transition at time $n+1$. $P_{\text{psinc}}^{n+1}$ is obtained by inserting $\pi^{n+1}_{\text{comp}}$ into the right-hand side of \eqref{eqn:P_conv_comp_to_psinc}. Moreover, there are two valid choices for $\pi$ in a pseudo-incompressible to compressible transition (Figure \ref{fig:timeline}): (1) $\pi_{\text{full}}$, \textit{i.e.}, $\pi$ obtained after the full $n$-to-$n+1$ time step; or (2) $\pi_{\text{half}}$, \textit{i.e.}, $\pi$ obtained after the $n+1$-to-\klein{$n+3/2$} half time-step.

$\pi_{\text{half}}$ is obtained from the solution of \eqref{eqn:half_time_rhs} with the solution of \eqref{eqn:half_time_lhs} as its input. The input to \eqref{eqn:half_time_lhs} are $\Psi^n$ and $(P\mathbf{v})^n$. This means that $\pi_{\text{half}}$ is recovered from the other quantities and is independent of $\pi$ at the previous time-level, so errors in the initialisation of $\pi$ are not propagated. By contrast, $\pi_{\text{full}}$ is obtained from the solution of \eqref{eqn:full_time}. The explicit \eqref{eqn:full_time_a} has $\pi$ as an input to the right-hand side $Q(\Psi^n;P^n)$. Therefore, $\pi_{\text{full}}$ propagates errors in the initialisation of $\pi$. Note that choice (2) entails solving an additional time-step in the pseudo-incompressible regime to obtain $\pi_{\text{half}}$ \revision{at the psuedo-incompressible to compressible blending time interfaces}.

\klein{In addition, choice (2) offers a conceptual advantage. The Exner pressure field in the pseudo-incompressible model is not controlled by an evolution equation but rather acts as a Lagrangian multiplier ensuring compliance of the velocity field with the divergence constraint \revision{at some fixed time}. Thus, a direct dependence of the pressure on its previous time level data, as occurs under option (1), is a numerical artefact that should be avoided.}

\subsection{Data assimilation and blending}
\label{subsec:dab}
A data assimilation engine is used to insert observations in the fully compressible configuration of the blended numerical framework. Prior to the assimilation procedure at time $t^n$, the forecast ensemble state vector $\{ \mathbf{x}^{f}_k \}^n$ for $k=1,\dots,K$ and a set of observations $\mathbf{y}_{\rm obs}^n$ are available. $K$ is the ensemble size. For vertical slice simulations with the \revision{full} compressible flow equations, the ensemble state vector is
\begin{equation}
    \{ \mathbf{x}^{f}_1, \dots \mathbf{x}^{f}_K \}^n = \{ \rho, \rho u, \rho w, P, \pi \}_{k=1,\dots,K}^n \in \mathbb{R}^{m \times K}.
    \label{eqn:state_vector}
\end{equation}
Here, the two-dimensional spatial grid has $(N_x \times N_z)$ cells and $m={(5 \times N_x \times N_z)}$. The observations of the momentum fields are
\begin{equation}
    \mathbf{y}_{\rm obs}^n = \{ (\rho u)_{\rm{obs}}, (\rho w)_{\rm{obs}} \}^n \in \mathbb{R}^l,
    \label{eqn:obs_vector}
\end{equation}
where the subscript \textit{obs} indicates that the data is obtained externally and is noisy and sparse, and $l = ({2 \times N_{\rm obs}(n)})$, with $N_{\rm obs}(n)$ the time-dependent dimension of the sparse observation space. The observation covariance $\pmb{\mathsf{R}}^n$ is determined by the observation noise.

The forward \klein{observation} operator $\mathcal{H}$ selects for each $\{ {\mathbf{x}}^{f}_k \}^n$ in \eqref{eqn:state_vector} the momenta $(\rho u)^n$ and $(\rho w)^n$ on the grid points corresponding to the sparse observations, thereby projecting $\mathbf{x}^{f, n}_k$ from the state space $\mathbb{R}^m$ into observation space $\mathbb{R}^l$, \textit{i.e.}
\begin{equation}
    \mathbf{y}^{f,n}_k = \mathcal{H}(\mathbf{x}^{f,n}_k) \in \mathbb{R}^l, \qquad  k=1,\dots,K.
    \label{eqn:state_vector_in_obs_space}
\end{equation}
The ensemble mean in observation space is computed as
\begin{equation}
    \bar{\mathbf{y}}^{f,n} = \frac{1}{K} \sum_{k=1}^K \mathbf{y}^{f,n}_k \in \mathbb{R}^{l}.
    \label{eqn:ens_mean_in_obs_space}
\end{equation}
A similar ensemble averaging is applied to obtain $\bar{\mathbf{x}}^{f,n}$. As observation localisation is used in the LETKF algorithm \citep{hunt2007}, only observations in a local region surrounding a given grid point are involved in its update. A localisation function is furthermore applied to the observations in the local region, \revision{see section \ref{subsec:da_setup} for more details on the setup}.

A Kalman gain $\pmb{\mathsf{K}}^n$ similar to \eqref{subeqn:ensKalman_gain} is obtained from the observation operator $\mathcal{H}$, the observation covariance $\pmb{\mathsf{R}}^n$, and an ensemble inflation factor $b$. As in the right-hand side of \eqref{subeqn:ensKalman_distance}, the distance of the forecast ensemble mean from the observations is computed with \eqref{eqn:obs_vector} and  \eqref{eqn:ens_mean_in_obs_space}. From these, a set of $K$ weight vectors $\{ \mathbf{w}^{a}_k \}^n$ is obtained, applied to \eqref{eqn:state_vector}, and added to $\bar{\mathbf{x}}^{f,n}$, updating the forecast ensemble to the analysis ensemble. Further details are given in Appendix A.

\begin{figure*}
\centering
\includegraphics[scale=1.0]{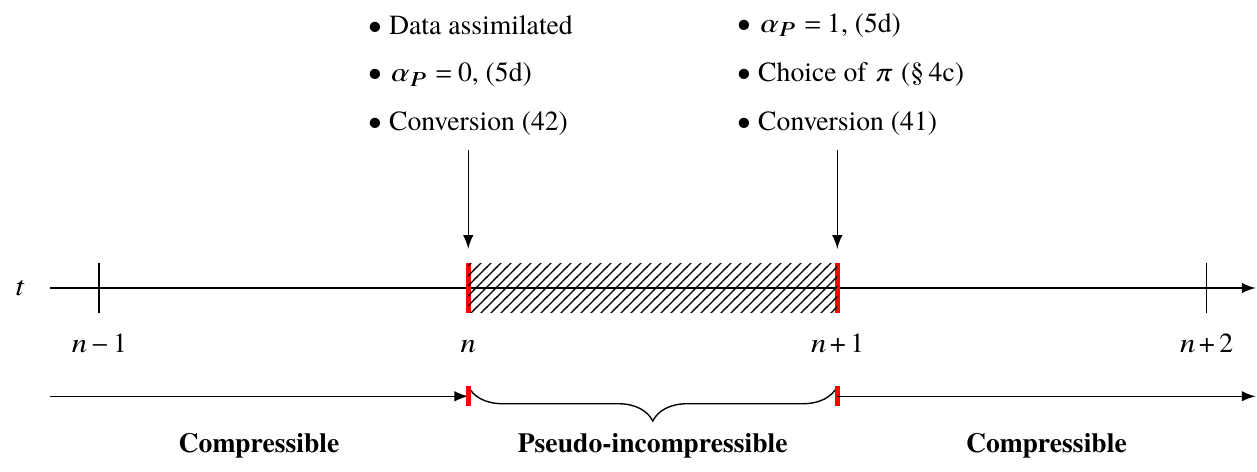}
\caption{Schematic of data assimilation with blending for data assimilated at time $t^n$. Blending time interfaces are in red, and the time-step spent in the pseudo-incompressible regime is shaded. See main text for full description. \revision{Numbers in brackets refer to equations, and (\S) denotes the section.}}
\label{fig:da_blending}
\end{figure*}

Once the assimilation procedure is completed, the model switches to the pseudo-incompressible limit regime and then back again to fully compressible until the next assimilation time. The process of switching back and forth between the model configurations exploits the blended numerical model to achieve balanced data assimilation and is termed \textit{blended data assimilation}.

\begin{figure}
\centering
\includegraphics[scale=1.0]{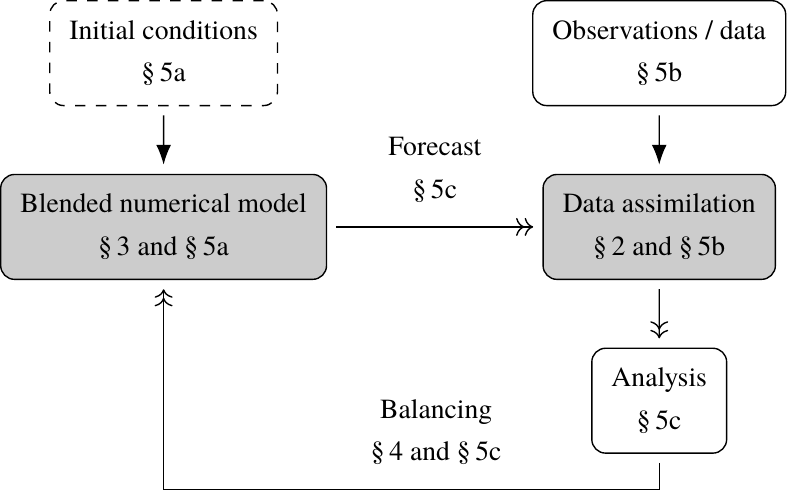}
\caption{Blended data assimilation workflow with the sections (\S) of this paper describing the algorithmic components. \revision{The initial condition (dashed outline) is used only once to start the simulation. For each assimilation window,  externally obtained observations / data are assimilated into the forecast and the algorithm loops through the components following the direction of the two-headed arrows.}}
\label{fig:overview_summary}
\end{figure}

In particular, if data \klein{are} assimilated into the compressible flow equations at time $n$, then compressible to pseudo-incompressible blending entails setting the switch $\alpha_P$ to $0$ and converting the quantity $P_{\text{comp}}$ with \eqref{eqn:P_conv_comp_to_psinc}. The solution is then \klein{propagated} in the pseudo-incompressible regime for a time-step, after which $\alpha_P$ is set back to 1, switching to the compressible flow equations. The quantity $P_{\text{psinc}}$ is reconverted by \eqref{eqn:P_conv_psinc_to_comp} using either $\pi_{\text{half}}$ or $\pi_{\text{full}}$. Figure \ref{fig:da_blending} summarises the procedure. Following the analysis from section \ref{subsec:pi_time_level}, the perturbation variable $\pi$ is reset after the half time-stepping in the solution of the full model, but not in the solution of the limit model. The blended data assimilation workflow is displayed in Figure \ref{fig:overview_summary}.

\revision{{As our principal strategy is to split measures of balancing the flow state from those of assimilating the data, we have not tuned the data assimilation procedures themselves in any way. Tuning its parameters may further improve balance, but as our balancing strategy is rather successful without this, the degrees of freedom of parameter tuning might be used more efficiently to achieve additional goals aside from the elimination of unphysical acoustic noise.}}

\section{Numerical results}
\label{sec:numerical_results}
The \revision{idealised} test cases of a travelling vortex and a rising warm air bubble are used to validate model performance in this section. To evaluate the effectiveness of the single time-step blended soundproof-compressible scheme, unbalanced states are initialised in the compressible flow equations for both test cases and the blended scheme is applied. The balance of the compressible solution with unbalanced initial states is evaluated by \revision{``probe measurements'', i.e., by time series of the flow variables at selected points in the domain,} and compared against \revision{analogous data extracted from the} soundproof solution \revision{ \citep{bok2014}}.

For blended ensemble data assimilation, an ensemble is generated \klein{by perturbing the initial conditions}\revision{. Then,} the blended scheme is applied after the assimilation of observations into the compressible flow equations\revision{ and repeated after each assimilation procedure}. The quality of balanced data assimilation is evaluated by root mean square errors with respect to a reference solution.

\subsection{Effectiveness of the improved blending strategy}
\label{subsec:effectiveness_blending}
\subsubsection{The travelling vortex experiment}
\label{subsubsec:vortex}

\begin{figure*}[ht]
    \centering
    \includegraphics[width=0.65\textwidth]{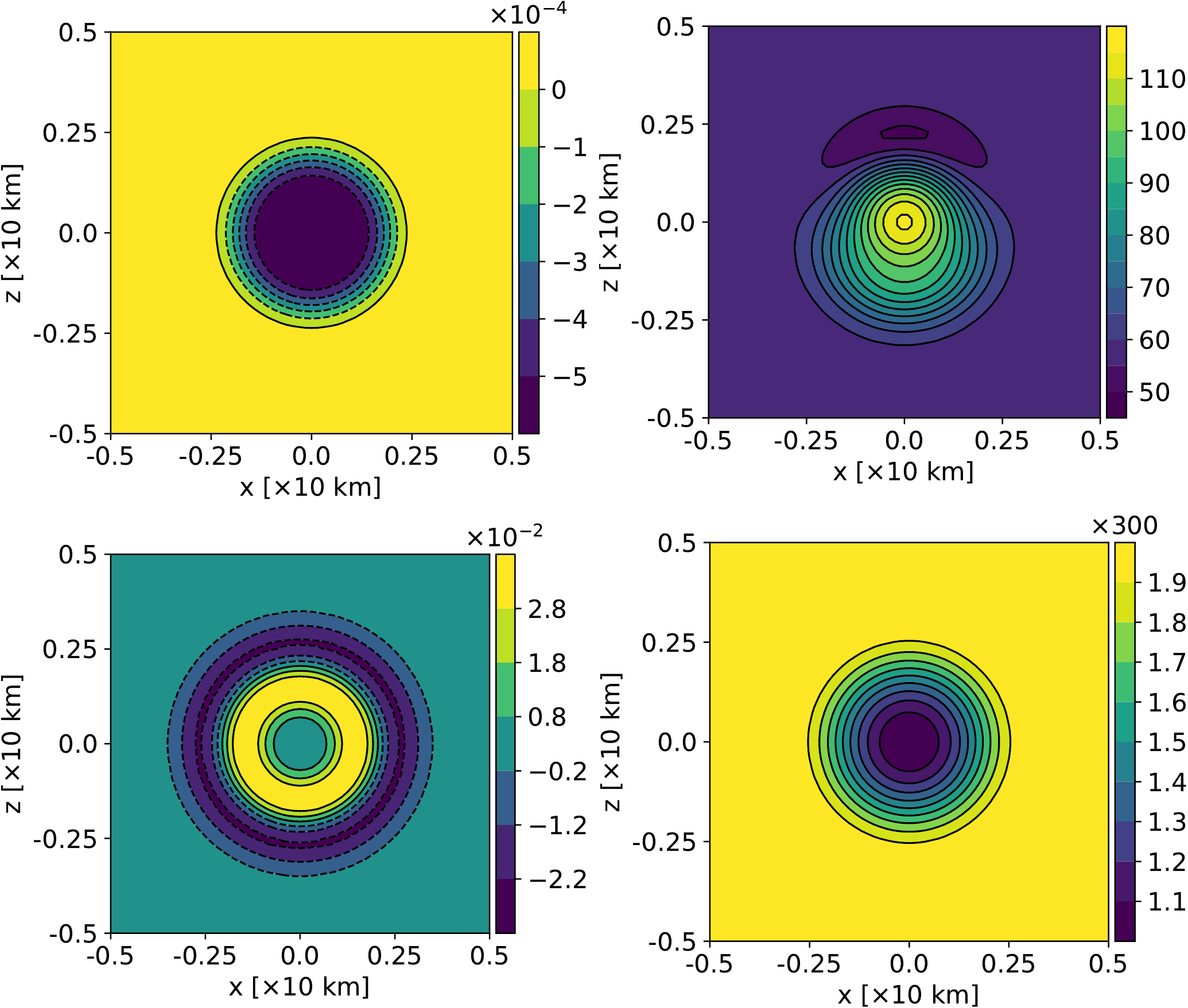}
    \caption{Travelling vortex initial balanced states: Exner pressure perturbation $\pi$; dimensionless contours in the range $[-5,0] \times 10^{-4}$ with interval of $10^{-4}$ (top left), horizontal momentum $\rho u$; contours in the range $[50, 115]$ kg m$^{-1}$ s$^{-1}$ with a $5$ kg m$^{-1}$ s$^{-1}$ interval (top right), vorticity; contours in the range $[-2.2,2.8] \times 10^{-2}$ s$^{-1}$ with a 1.0 $\times 10^{-2}$ s$^{-1}$ interval (bottom left), and potential temperature $\Theta$; contours in the range $[1.1,1.9] \times 300$ K with a 0.1 $\times  ~300$ K interval (bottom right). Negative contours dashed.}
    \label{fig:initial_vortices}
\end{figure*}

The travelling vortex test case of \cite{kadioglu2008} with $f=0.0$ s$^{-1}$ and $g=0.0$ m s$^{-2}$ is considered in the domain $x=[-5.0\text{ km},5.0\text{ km}]$, $z=[-5.0\text{ km},5.0\text{ km}]$ with doubly periodic boundary conditions and a background wind with velocity $100$ m s$^{-1}$ in both directions (Figure \ref{fig:initial_vortices}). The time-step size is constrained by advective $\mathrm{CFL} = \| \mathbf{u} \| \Delta t / \Delta x = 0.45$ on a $(64 \times 64)$ grid. The choice of reference units yield\revision{s} $\mathrm{Ma} \approx 0.341$. \klein{Note, however, that while the background wind Mach number is relatively large, the superimposed vortex has a maximum flow velocity of $25\text{\ m s$^{-1}$}$ with $\mathrm{Ma}_{\text{vort}} = 0.076$, so that the low Mach number analysis of section~\ref{ssec:Conversion} is justified.}

In order to gauge the performance of the improved blended model, probe measurements of the \revision{full} pressure increments \revision{$\delta p$} are taken, defined, e.g. at time-level $n$, as
\begin{equation}
    \revision{\delta p^{n} = p^{n+1} - p^{n}.}
    \label{eqn:pressure_increment}
\end{equation}
at the center (0.0 km, 0.0 km). The first increment \revision{$\delta p^{0}$} corresponds to a spinup adjustment and is therefore omitted in the plots, as done in \cite{bok2014}. 

The distance in the pressure perturbation result of a \revision{given} run compared to the reference pseudo-incompressible run is quantified by the relative error $E_\nu$,
\begin{equation}
    \revision{E_{\nu} = \frac{\lVert \delta p_{\nu} - \delta p_{\rm psinc} \rVert_2}{\lVert \delta p_{\rm psinc} \rVert_2},}
    \label{eqn:deltas}
\end{equation}
where $\nu = \rm{b}$ for the blended run and $\nu = \rm{c}$ for the imbalanced compressible run. $\lVert \, \cdot \, \rVert_2$ is the 2-norm taken over the probe measured time series of \revision{$\delta p$}. 

\revision{An imbalanced initial} state is created by setting $P=10^{5}$ Pa and $\pi=0.0$ over the whole domain for the \revision{full} compressible flow equations \revision{\eqref{eqn:governing}} with $\alpha_P = 1$. \revision{This state is propagated for} one time-step in the limit pseudo-incompressible regime followed by the rest of the time-steps in the fully compressible model. The blending scheme in section \ref{sec:blending_scheme} is used to transition between the model regimes.

\begin{figure*}[!ht]
    \centering
    \includegraphics[width=0.65\textwidth]{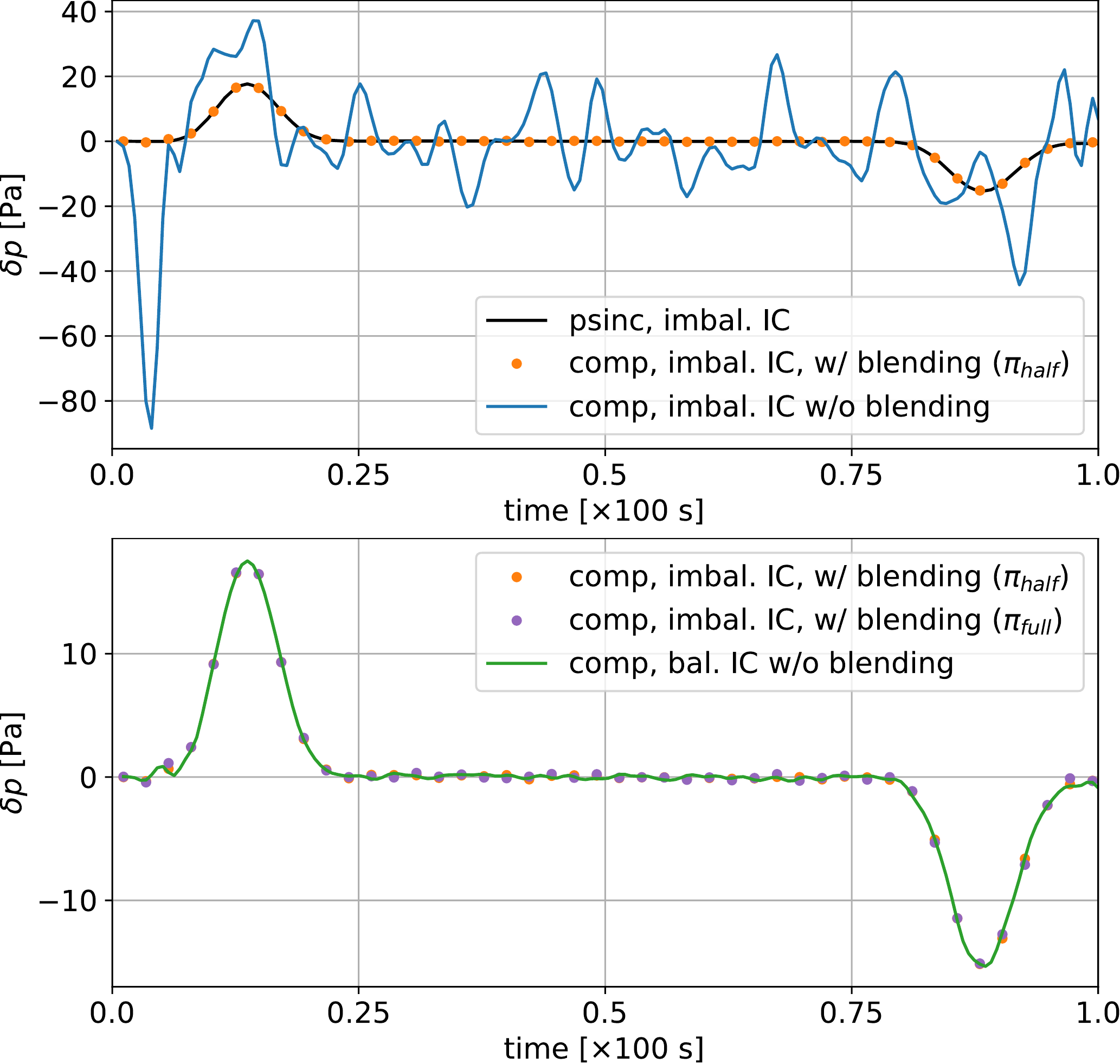}
    \caption{Travelling vortex: effect of blending for imbalanced initial states \revision{on the time series of temporal increments of the full pressure $\delta p$ at location $(x,z) = (0 \text{ km},0 \text{ km})$.} Top: comparison between a blended run using  $\pi_{\text{half}}$ (orange), a run without blending (blue), and the reference solution from the pseudo-incompressible model (black). Bottom: comparison of blended runs using  $\pi_{\text{half}}$ (orange) and $\pi_{\text{full}}$ (purple), and the compressible solution with balanced initial states (green). The blended runs are with one time-step spent in the pseudo-incompressible regime.}
    \label{fig:comparison_w_imbal}
\end{figure*}

For these imbalanced initial states, a compressible run with blending is compared \klein{with} \revision{a compressible run} without blending and with a pseudo-incompressible run (top panel in Figure \ref{fig:comparison_w_imbal}). Fast acoustic modes are filtered from the blended solution and the result is indistinguishable from the limit pseudo-incompressible reference solution, save for an initial adjustment in the first time-step. Blending is able to recover the dynamics of the balanced state.

A close-up (bottom panel of Figure \ref{fig:comparison_w_imbal}) compares the blended runs with choices of $\pi_{\text{half}}$ and $\pi_{\text{full}}$ from section \ref{subsec:pi_choice} against a run with \klein{the balanced initial state obtained from the known exact compressible vortex solution}. The blended runs are as good as, and the $\pi_{\text{half}}$ run slightly closer to, the balanced compressible run. The relative error of the blended run with respect to the reference balanced run is \revision{0.0285} using $\pi_{\text{half}}$ and \revision{0.0393} using $\pi_{\text{full}}$. This corroborates the insight from section \ref{subsec:pi_choice} that $\pi_{\text{half}}$ is a better choice. The choice of $\pi_{\rm half}$ is used from here on.

\subsubsection{The rising bubble experiment}
\label{subsubsec:bubble}

\begin{figure*}
    \centering
    \includegraphics[width=0.65\textwidth]{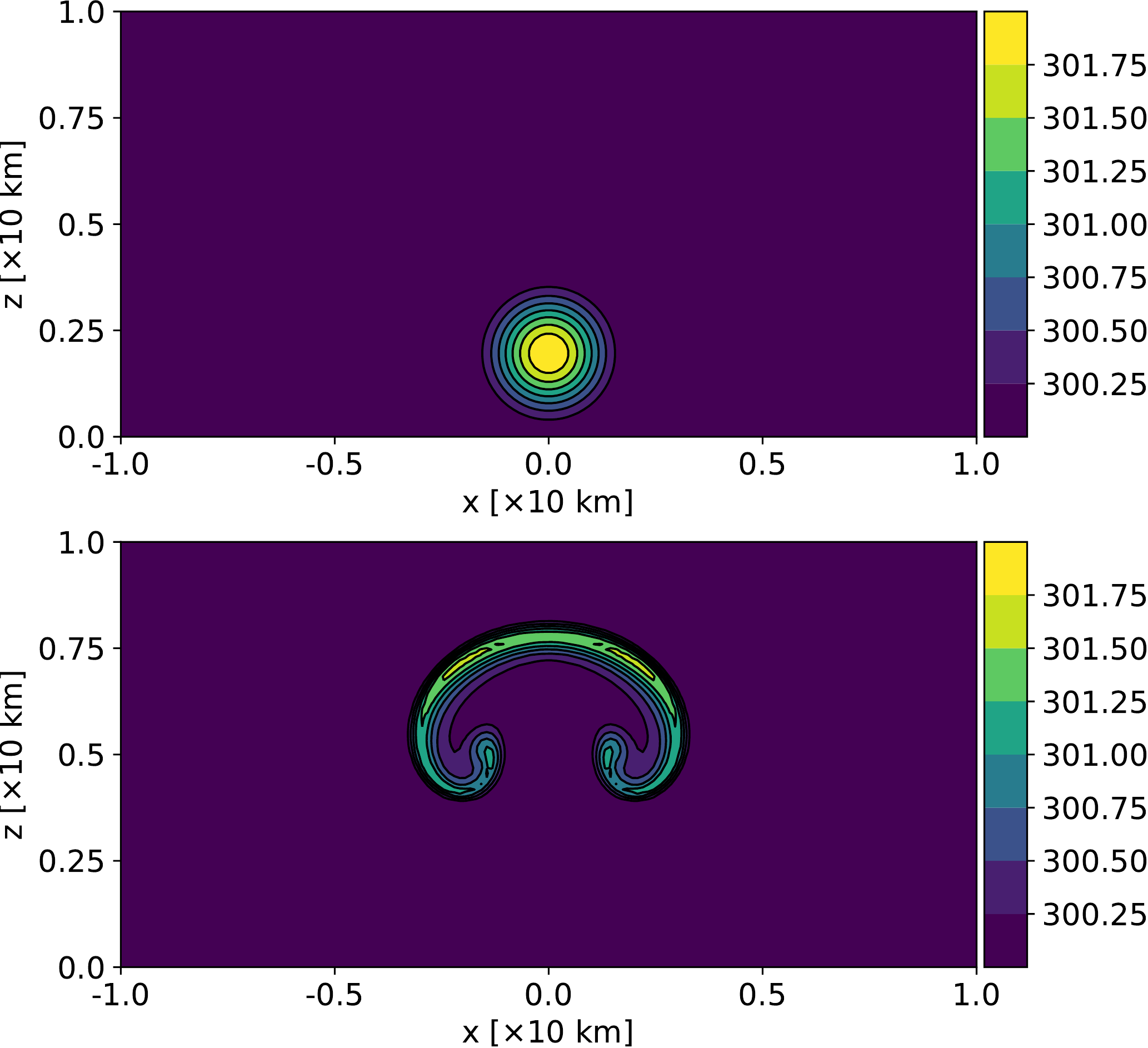}
    \caption{Rising bubble: potential temperature at initial time, $t=0$ s \revision{(top)} and final time $t_{\text{fin}}=1000$ s \revision{(bottom)}; contours in the range $[300.25,301.75]$ K with a $0.25$ K interval.}
    \label{fig:initial_rb}
\end{figure*}

The second test consists of a gravity-driven thermal flow with $f=0.0$ s$^{-1}$ initialised as a bubble-shaped positive potential temperature perturbation $\delta \Theta$, on a constant isentropic background with $\Theta_0 = 300$ K in a $[-10.0\text{ km},10.0\text{ km}]\times[0.0\text{ km},10.0\text{ km}]$ domain, with periodic boundaries in $x$ and no-flux in $z$ \revision{\citep{MendezCaroll1994,klein2009asymptotics,bok2014}}. The dimensionless perturbation $\delta \Theta$ is defined by
\begin{equation}\label{eq:InitialThermal}
    \delta \Theta = \frac{2\, \text{K}}{\Theta_0} \cos \left( \frac{\pi }{2} r \right),
\end{equation}
where
\begin{equation}
    r = \frac{1}{r_0} \sqrt{x^2 + (z-z_0)^2},
\end{equation}
$r_0=2.0$ km is the initial radius of the bubble, and $z_0=2.0$ km the initial vertical displacement of the bubble. The choice of reference units yields $\mathrm{Ma} \approx 0.0341$. The models are run on a grid with $(160 \times 80)$ cells to a final simulation time of $1000.0$ s.

The initial pressure fields are \klein{set to reflect a horizontally homogeneous hydrostatic pressure field $\overline{p}(z)$ based on $\Theta_0$ and initial condition $\overline{p}(0) = 10^5\, \text{Pa}$, with $\revision{\pi^\prime} \equiv0.0$. These pressure data are imbalanced, however, with respect to the perturbed initial potential temperature $\Theta_0 + \delta \Theta$, see \eqref{eq:InitialThermal}}. Potential temperature at the initial and final time are depicted in Figure \ref{fig:initial_rb}.

\begin{figure*}[!ht]
    \centering
    \makebox[\textwidth][c]{
    \includegraphics[width=1.1\textwidth]{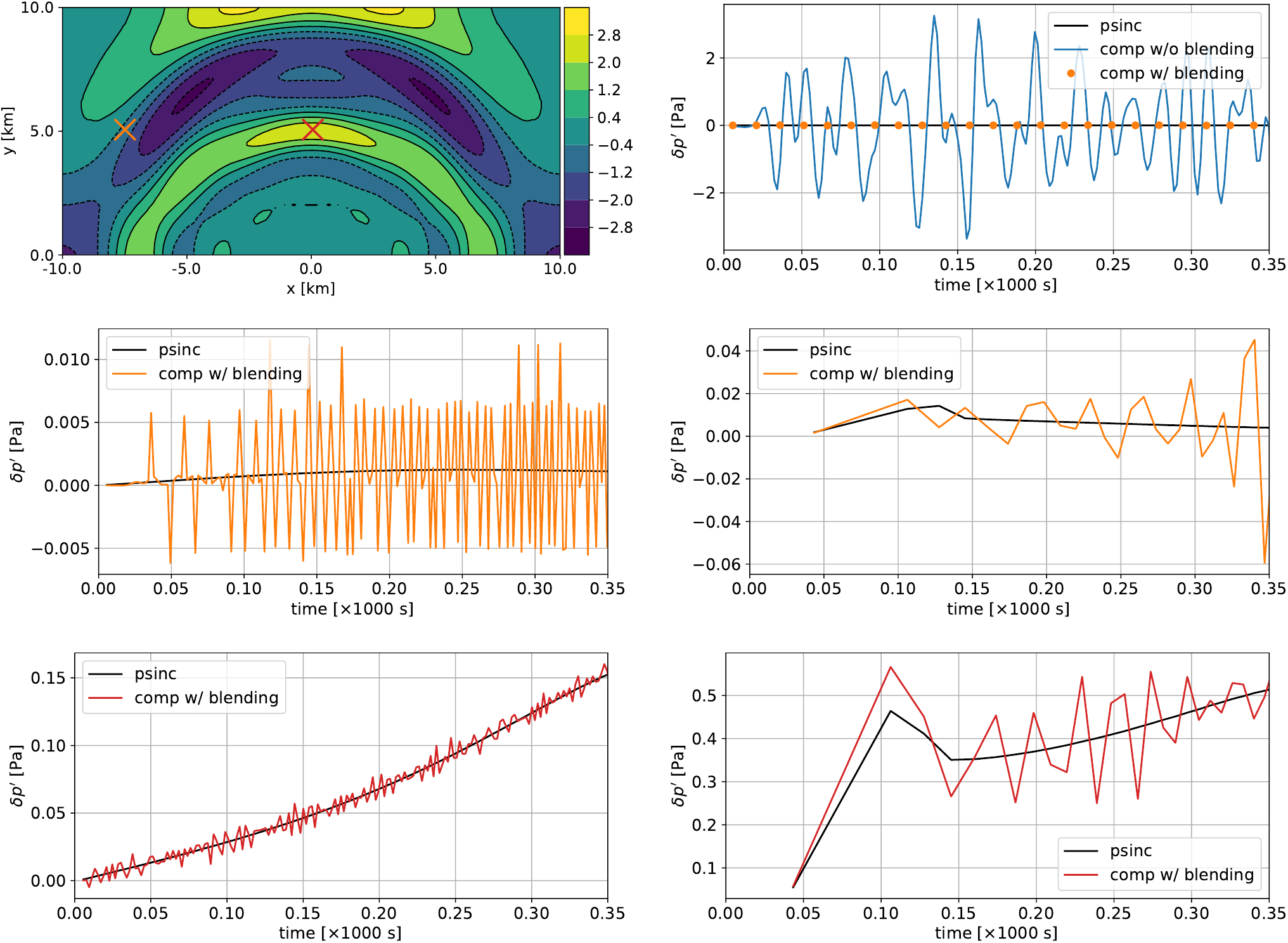}
    }
    \caption{Time increments $\delta p^\prime$ of pressure perturbation in the rising bubble experiment. Top left: $\delta p^\prime$ value at time-step 14 ($t=26.6$ s) for the compressible model; contours in the range $[-2.8,2.8]$ Pa with an interval of $0.8$ Pa, negative contours dashed. The orange cross marks $(x,z)=(-7.5,5)$ km and the red cross $(x,z)=(0,5)$ km. Top right: history of $\delta p^\prime$ over the first $350$ s measured at $(x,z)=(-7.5,5)$ km for the compressible model (blue), the pseudo-incompressible model (black) and the blended model with one pseudo-incompressible time-step (orange).  Middle and bottom panels: history of $\delta p^\prime$ over the first $350$ s measured at $(x,z)=(-7.5,5)$ km (middle) and at $(x,z)=(0,5)$ km (bottom). Pseudo-incompressible solution (black) and blended solution with one time-step spent in the pseudo-incompressible regime (orange or red corresponding to the probe marker in the top left panel).
    The top, middle left, and bottom left panels use a constant time-step $\Delta t = 1.9$ s. The middle right and bottom right panels use $\Delta t= 21.69$ s for the first two time-steps and then a $\Delta t$ determined by advective $\mathrm{CFL}=0.5$.
    }
    \label{fig:rising_bubble}
\end{figure*}

The initial stages of the bubble evolution are compared for the compressible, pseudo-incompressible and \revision{one-step} blended runs in Figure \ref{fig:rising_bubble}.  As the initial state is not hydrostatically balanced, pressure waves \klein{propagate} in the compressible configuration (top left panel) as seen in a time series of pressure perturbation increment at $(x,z)=(-7.5,5)$ km (orange cross in the top left panel and \revision{blue line in the} top right panel). The acoustics are absent in the soundproof configuration (top right panel, black line) and in the single time-step blended soundproof-compressible configuration (orange dots).

Next, the blended run and the pseudo-incompressible run are compared in more detail (Figure \ref{fig:rising_bubble}, middle and bottom panels) with pressure perturbation increment probe measurements \revision{$\delta p^\prime$, and $p^\prime = p-\bar{p}$. The probes are located} at $(x,z)=(-7.5,5)$ km (middle panels) and at $(x,z)=(0,5)$ km (red cross in the top left panel and bottom panels in Figure \ref{fig:rising_bubble}), both with a constant small time-step $\Delta t=1.9$ s (top, middle left and bottom left panels) and for larger, advective CFL-constrained time-steps (middle right and bottom right panels,  $\mathrm{CFL}=0.5$ and $\Delta t = 21.69$ s for the first two time-steps). Away from the bubble trajectory (middle panels), the pressure perturbation increment due to the rising bubble and the remnants of the background acoustics from blending are comparable in amplitude. Larger amplitudes are observed with the blended model and the larger time step (middle right panel), but they are still very small compared to the fully compressible run (note the different range on the vertical axes between the top right and middle right panels). On the bubble trajectory (bottom panels), the pressure perturbation increment due to the rising bubble dominates and the solutions are almost identical.

\begin{table}
\caption{\footnotesize{Errors $E_{\rm c}$ and $E_{\rm b}$ (see text for definitions) of the time series of $\delta p^\prime$ in $[0,\,1000]$ s relative to the reference pseudo-incompressible run (middle and bottom panels of Figure \ref{fig:rising_bubble}). $\Delta t_{AC}=1.9$ s, $\Delta t_{ADV}$ is determined by advective $\rm{CFL}=0.5$ and $\Delta t_{ADV} = 21.69$ s for the first two time-steps. Probe location $(-7.5,5)$ km corresponds to the orange marker and orange lines in Figure \ref{fig:rising_bubble} and $(0,5)$ km to the red markers and red lines.}}\vspace{3mm}
        \centering\footnotesize
    \begin{tabular}{ccccc}
    \toprule\midrule
    probe location               & $\Delta t$ & $E_{\rm c}$ & $E_{\rm b}$ & $E_{\rm c}$/$E_{\rm b}$ \\ \midrule
    \multirow{2}{*}{$(-7.5,5)$ km} & $\Delta t_{AC}$      & 413.1822    & 1.4820 & 278.80      \\ \cline{2-5} 
        & $\Delta t_{ADV}$      & 109.7538    & 3.9034 & 28.12      \\ \midrule
    \multirow{2}{*}{$(0,5)$ km}    & $\Delta t_{AC}$      & 10.1804     & 0.0311 & 327.34     \\ \cline{2-5} 
        & $\Delta t_{ADV}$      & 2.8231      & 0.1016 & 27.79      \\ \bottomrule
    \end{tabular}
\label{tab:error}
\end{table}

Throughout the runs, a single time-step spent in the soundproof pseudo-incompressible regime largely filters out the fast acoustic imbalances of the compressible run (not shown in the middle and bottom panels of Figure \ref{fig:rising_bubble}). This is quantified by comparing the relative errors with respect to the reference pseudo-incompressible run for the compressible run, $E_{\rm{c}}$, and for the blended run, $E_{\rm{b}}$, defined in \eqref{eqn:deltas} and shown in Table \ref{tab:error}. $E_{\rm{b}}$ is more than 25 times smaller than $E_{\rm{c}}$ for the large time-step case, and more than two orders of magnitude smaller for the small time-step case.

\revision{We also remark that a probe measurement of the \textit{full pressure} time increment $\delta p$ differs slightly between the reference pseudo-incompressible run and the one-step blended run (not shown). The difference is due to the time-dependence of the hydrostatically balanced background pressure $\bar{p}$ in the blended run. However, the computed values of the \textit{pressure perturbation} time increment $\delta p^\prime$ are remarkably similar in the two runs (black line and orange dots in top right panel of Figure \ref{fig:rising_bubble}). We can thus conclude that blending recovers balanced dynamics irrespective of small compressibility-induced variations of the background pressure $\bar{p}$.

In view of these results, b}lending can be employed as an effective means to achieve the balanced initialisation of data within a fully compressible model. The single time-step balancing capability in the model presented here substantially improves on the performance of \cite{klein2014using} and \cite{bok2014}, whose blended models achieved smaller reductions in amplitude compared to the fully compressible case and needed several time steps in the limit regime.

\subsection{Ensemble data assimilation and blending: setup}
\label{subsec:da_setup}
\subsubsection{Travelling vortex setup}
To combine blending with data assimilation as described in Section 4d, an ensemble is generated by perturbing the initial vortex center position $(x_c,z_c)$ within the open half interval of $[-1.0 \text{ km}, 1.0 \text{ km})$ for both $x_c$ and $z_c$. The vortex is then generated around this center position such that the full vortex structure is translated. Ten such samples are drawn and they constitute the ensemble members. An additional sample is drawn and solved with the full model for the balanced initial condition. This run, denoted by \textit{obs}, is used to generate the artificial observations. Another run identical to this additional obs sample is made. \klein{This time, \revision{however,} blending for the first time-step is applied and \revision{this run} is \klein{considered} the \textit{truth} in the sequel}. This is to correct for any errors in the initialisation of $\pi$, as discussed in section \ref{subsec:pi_choice}.

This choice of generating the truth and obs through a perturbation of the initial condition is such that the ensemble mean does not coincide with the truth. Otherwise, ensemble deflation alone \revision{would be} sufficient to make the ensemble converge towards the truth, see also \cite{lang2017data}.

The observations are taken from the obs run every 25 s -- only a tenth of the grid points are observed and these are drawn randomly. \revision{Sparse observation grid points are randomly drawn as follows: A Boolean mask selecting for a tenth of the grid points is generated where if necessary, a ceiling function is applied to obtain an integer number of grid points selected. The entries of the mask are then shuffled using the algorithm by  \citet{fisher1953statistical} and the Boolean mask is applied to the obs array to obtain the sparse observations.} This deviates from a more realistic situation where observations and grid points do not coincide. To simulate measurement noise, Gaussian noise with standard deviation equal to 5\% of the peak-to-peak amplitude of the obs quantity at the given time is added independently to each of the observed grid points. A similar method of generating artificial observations was used in, for example, \cite{bocquet2011ensemble, harlim2005local} for the Lorenz-63 and Lorenz-96 models.

The regions \klein{for localised data assimilation} are of size $(11 \times 11)$ grid points and only observations within \klein{such a patch} are considered for analysis \klein{operations at the respective central grid point}. A localisation function corresponding to a truncated Gaussian function is applied such that observations farther from the grid point under analysis have less influence, and that the influence decays smoothly \klein{towards the edges of the localisation subdomain, \revision{where it is abruptly truncated to zero}}. \klein{No ensemble inflation is applied in this case}.

\begin{figure}[!ht]
    \centering
    \includegraphics[width=0.65\textwidth]{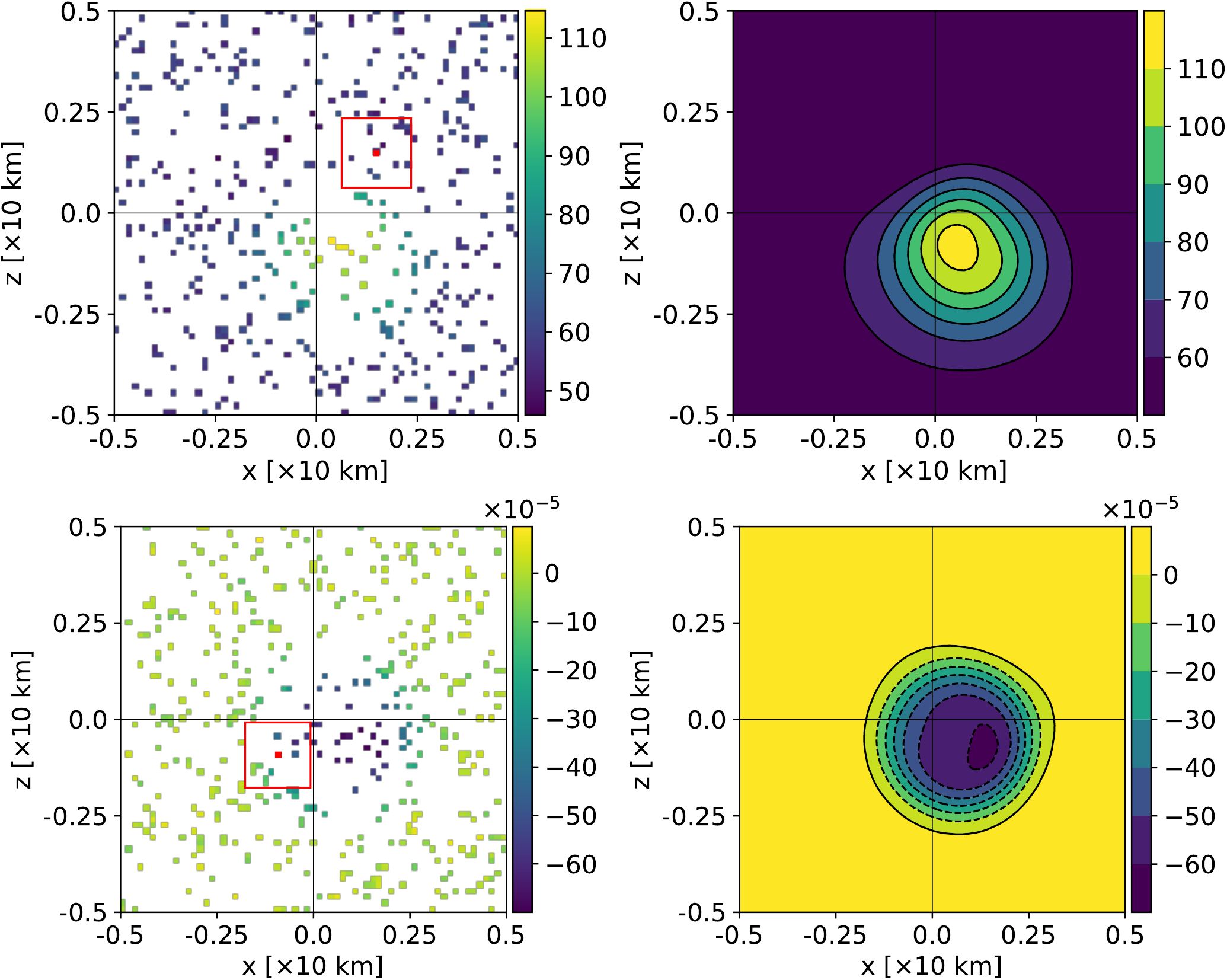}
    \caption{Travelling vortex: sparse noisy observations (left panels) and truths (right panels) at $t=300.0$ s. Top: horizontal momentum $\rho u$; contours in range $[60,110]$  kg m$^{-1}$ s$^{-1}$ with a $10$ kg m$^{-1}$ s$^{-1}$ interval. Bottom: Exner pressure perturbation $\pi$; dimensionless contours in range $[-60,0] \times 10^{-5}$ with an interval of $10^{-4}$. Negative contours are dashed. The red squares illustrate, for an example grid point (in red), the observations considered in the local $(11 \times 11)$ grid points region.}
    \label{fig:obs_truth}
\end{figure}

Examples of the observations and truths used in the generation and evaluation of the experiments with data assimilation are displayed in Figure \ref{fig:obs_truth}. \klein{Notice that we run one test with observations of the momentum fields only, and another test with observations of the full set of variables.}

The ten ensemble members \klein{in each of these tests} are initialised with balanced states, and blending is applied for the first time-step when the model runs in the pseudo-incompressible configuration. The ensemble is then solved forward in time with the fully compressible model. Data from the generated observations are assimilated every 25~s. The immediate time-step after the assimilation procedure is solved in the pseudo-incompressible limit regime while the rest of the time-steps in the assimilation window are solved using the fully compressible model. Conversions according to the blending scheme in section~\ref{sec:blending_scheme} are \klein{employed when switching} back and forth between the full and limit models. Furthermore, the choice of $\pi_{\text{half}}$ is used \revision{(cf.~the discussion in section~\ref{subsec:pi_choice}.)}. The ensemble solved with both data assimilation and blending is abbreviated as EnDAB.

The setup is repeated for two \klein{additional ensembles and each observation scenario}, one where data are still assimilated but no blending is performed (EnDA), and another where neither data assimilation nor blending are performed (EnNoDA). EnNoDA and EnDA constitute an identical twin experiment \citep{reich2015probabilistic, lang2017data}, through which the effects of data assimilation can be evaluated. EnDA along with EnDAB constitute yet another identical twin experiment, which evaluates the performance of blending.

\subsubsection{Rising bubble setup}
The rising bubble ensemble spread is generated by randomly \klein{modifying} the \klein{maximum of the} potential temperature perturbation $\delta \Theta$ in the open half interval $[2.0 \text{ K}, 12.0 \text{ K})$. The ensemble comprises ten members. While the \klein{relative spread of the} temperature perturbation is large \klein{with this setup}, the ensemble spread of the bubble position at the final time of the simulation, $t_{\text{fin}}=1000.0$ s, is only moderate.

\begin{figure}[!ht]
    \centering
    \includegraphics[width=0.65\textwidth]{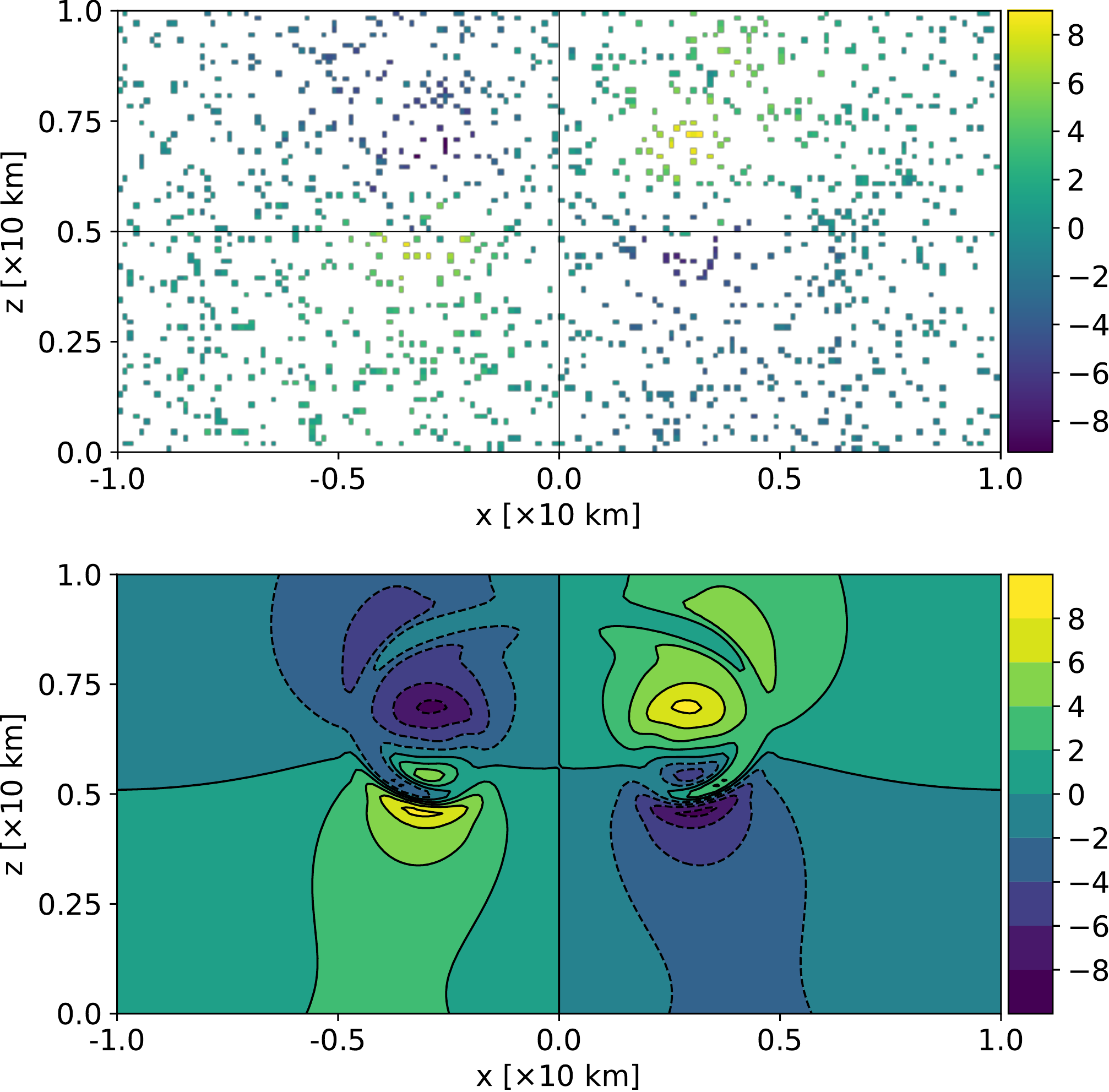}
    \caption{Rising bubble: horizontal momentum field $\rho u$ at $t=1000.0$ s. Sparse and noisy observations (top) and truth (bottom); contours in range $[-8,8]$ kg m$^{-1}$ s$^{-1}$ with a $2$ kg m$^{-1}$ s$^{-1}$ interval, negative contours dashed.}
    \label{fig:rb_obs_truth}
\end{figure}

An additional sample is drawn for the obs and the truth, which are identical in this setup. Blending is applied to the first time-step of the obs and the truth, obtaining a balanced solution. \klein{As the rising bubble flow fields evolve rather slowly in the beginning, data are only assimilated from $t=500.0$~s onwards. Observations of the momentum field are then assimilated every $50.0$ s}. As with the vortex experiments, only a tenth of the grid points are observed, noise with standard deviation 5\% of the peak-to-peak amplitude is added, and localisation within an $(11 \times 11)$ grid points region is applied. A localisation function corresponding to the truncated the Gaussian function is applied and the ensemble is not inflated. Examples of the observation and truth are given in Figure \ref{fig:rb_obs_truth}. Three ensembles corresponding to the EnNoDA, EnDA, and EnDAB settings, with 10 members each, are generated, \klein{but only one set of experiments involving assimilation of the momentum field only is pursued}.

Note that as the ensembles and the observations are generated with balanced initial conditions, any noise present in the simulation results is the result of the data assimilation procedure. Table \ref{tab:experimental_setup} summarises the details of the data assimilation-related experimental setup for both test cases. 

\begin{table}[!ht]
\caption{\footnotesize{Assimilation-related experimental parameters. $K$ is the ensemble size, $b$ the ensemble inflation factor, $t_{\rm first}$ the first assimilation time, $\Delta t_{\rm obs}$ the observation interval, $\psi_{\rm assimilated}$ the set of quantities assimilated, $(N \times N)_{\rm local}$ the size of the local region, $f_{\rm local}$ the type of localisation function, $\eta_{\rm obs}$ the observation noise, $\sigma$ the standard deviation of the Gaussian noise, $A$ the peak-to-peak amplitude of the quantity observed, obs$_{\rm sparse}$ the sparsity of the observations, $\Delta t_{\rm blending}$ the number of initial time-steps spent in the limit model regime. $\pi$ choice, used in the initialisation of $\Delta t_{\text{blending}}$, is either $\pi_{\rm half}$ or $\pi_{\rm full}$, more details in section \ref{subsec:pi_choice}.}}\vspace{6 mm}
\centering
\small
\begin{tabularx}{1.0\textwidth}{@{\extracolsep{\fill}}llcccc@{\extracolsep{\fill}}}
\toprule\midrule
\textbf{Test case}  &  & & & \textbf{Vortex} & \textbf{Bubble}  \\ \midrule
\multirow{2}{*}{\textbf{Ensemble}} & & $K$ 
&    & \multicolumn{2}{c}{10 members} \\
 & & $b$
  &  & \multicolumn{2}{c}{1.0} \\ \midrule
\multirow{7}{*}{\textbf{Observations}} & &
    $t_{\rm first}$ [s]
&    & $25.0$
    & $500.0$ \\
 &  & $\Delta t_{\rm obs}$ [s]
    & & $25.0$ 
    & $50.0$ \\
 & &   $\psi_{\rm assimilated}$
  &  & \stackanchor{$\{ \rho u, \rho w \}$ or}{$\{ \rho, \rho u, \rho w, P, \pi \}$}
    & $\{ \rho u, \rho w \}$ \\
  & & $(N \times N)_{\rm local}$ & & \multicolumn{2}{c}{$(11\times11)$ grid points} \\
& &   $f_{\rm local}$
 &   & \multicolumn{2}{c}{Truncated Gaussian} \\
 & &  $\eta_{\rm obs}$
   & & \multicolumn{2}{c}{Gaussian with $\sigma=0.05A$ } \\ 
& &   obs$_{\rm sparse}$
 &   & \multicolumn{2}{c}{One in 10 grid points} \\ \midrule
 \multirow{2}{*}{\textbf{Blending}} 
 & &  $\Delta t_{\rm blending}$
  &  & \multicolumn{2}{c}{A single blended time-step} \\
  & & $\pi$ choice
&    & \multicolumn{2}{c}{$\pi_{\text{half}}$}\\
    \bottomrule
\end{tabularx}
\label{tab:experimental_setup}
\end{table}

\subsubsection{Evaluation of data assimilation}
The quality of data assimilation is evaluated by a spatially and ensemble averaged root mean square error (RMSE) from the truth. This is given by
\begin{equation}
\footnotesize{
\text{RMSE}(\psi) = \sqrt{ \frac{1}{K} \frac{1}{N_x \times N_z} \sum_k^K \sum_{i,j}^{N_x,N_z} \left[ \psi^{\text{ensemble}}_k(x_i,z_j)  - \psi^{\text{truth}}(x_i,z_j) \right]^2 },
}
\label{eqn:rmse}
\end{equation}
where $k=1,\dots,K$ indexes the ensemble members and $i=1,\dots,N_x$ and $j=1,\dots,N_z$ the number of grid points in the $(x,z)$ coordinates. $\psi$ is the set of quantities $\{ \rho, \rho u, \rho w, P, \pi \}$.

\subsection{Ensemble data assimilation and blending: results}
\label{subsec:da_results}
\subsubsection{Travelling vortex}

\begin{figure*}[!ht]
    \centering
    \makebox[\textwidth][c]{
    \includegraphics[width=1.25\textwidth]{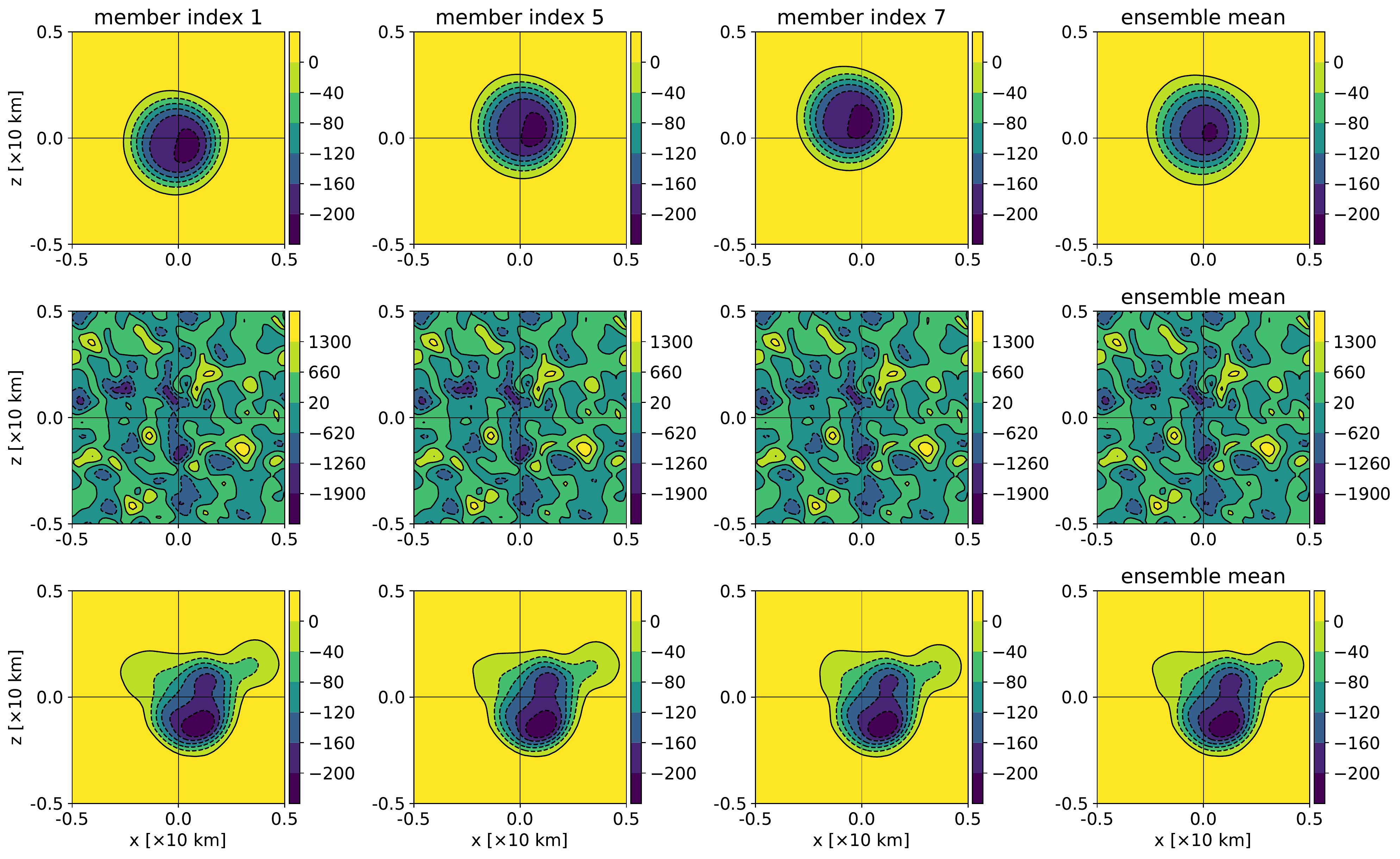}
    }
    \caption{Travelling vortex: snapshots of pressure perturbation $p^\prime$. Ensemble members with index 1 (first column), 5 (second column), 7 (third column) and ensemble mean (fourth column) at $t=300.0$ s with all quantities $\{\rho, \rho u, \rho w, P, \pi\}$ assimilated. Top row: EnNoDA run; contours in range $[-200, 0]$ Pa with a $40$ Pa interval. Middle row: EnDA run; contours in range $[-1900,1300]$ Pa with a $640$ Pa interval. Bottom row: EnDAB run; contours in range $[-200,0]$ Pa with a $40$ Pa interval. Negative contours are dashed.}
    \label{fig:euler_enses_3.0}
\end{figure*}

Figure \ref{fig:euler_enses_3.0} depicts the ensemble snapshots for the vortex case with all quantities observed and assimilated. \revision{EnNoDA acts as the control ensemble, and the top row depicts its solutions for the travelling vortex without data assimilation and blending}. While the center position of the vortex for each ensemble member is perturbed, the ensemble mean vortex position (fourth column) is centered around the origin. This is in line with the conditions used to generate the initial ensemble. With data assimilation, EnDA (middle row), the balance is lost and at final time the vortex structure is not preserved. Data assimilation and blending, EnDAB (bottom row), recovers the balanced solution and the vortex structure is preserved after three periods of revolution. Moreover, comparing with Figure \ref{fig:obs_truth}, the effect of data assimilation becomes obvious. The center position of the EnDAB ensemble mean is in the lower right quadrant, closer to that of the observation and the truth.

\begin{figure*}[ht]
    \centering
    \makebox[\textwidth][c]{
    \includegraphics[width=1.25\textwidth]{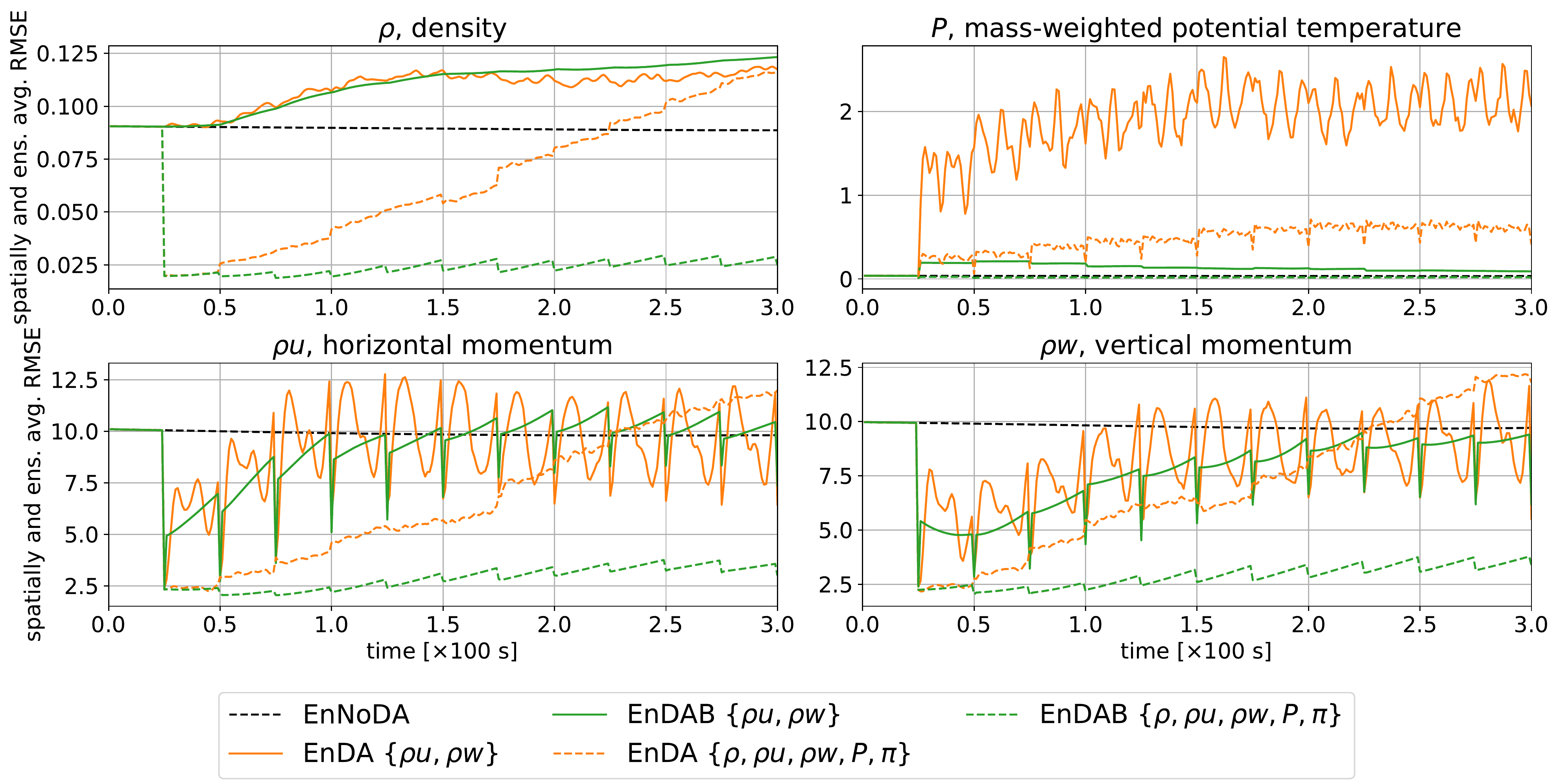}
    }
    \caption{Travelling vortex: EnNoDA run (black), EnDA run (orange), EnDAB run (green). Assimilated quantities are $\rho u$ and $\rho w$ (solid lines) and $\rho, \rho u, \rho w, P, \pi$ (dashed lines). Spatially and ensemble averaged RMSE from $t=0.0$ s to $300.0$ s for density $\rho$ (top left, [kg m$^{-2}$]), mass-weighted potential temperature $P$ (top right, [kPa]), and momenta $\rho u$, $\rho w$ (bottom left and right, [kg m$^{-1}$ s$^{-1}$]). The RMSE of the initial ensemble is omitted.}
    \label{fig:euler_rmse_all}
\end{figure*}

\revision{Data assimilation without blending (EnDA, orange lines in Figure \ref{fig:euler_rmse_all})}
leads to a jump in the RMSE in the thermodynamic $P$ variable upon the first assimilation at $t=25$ s. After that, the error stays relatively constant. Appendix B shows that the error jump quantifies the imbalance introduced by the data assimilation procedure.

\revision{Assimilating the momentum fields alone is insufficient and the RMSE in the solution (solid lines in Figure \ref{fig:euler_rmse_all}) is larger than in the reference EnNoDA run. As expected, EnDAB provides a smoother solution over time as the error does not oscillate\revision{, yet} ensemble spread and RMSE are comparable in these runs (not shown).} This test includes a strong axisymmetric potential temperature variation (Figure~\ref{fig:initial_vortices}), and the potential temperature is an advected quantity not corrected by momentum data assimilation. Therefore, the initially tight correlation of the velocity and potential temperature variations gets destroyed in the course of data assimilation. Since the potential temperature is fluid dynamically active through the generation of baroclinic torque, the flow fields of the ensemble members increasingly deviate from their reference as a consequence.

\revision{Assimilating all the quantities yields an improvement (dashed lines in Figure \ref{fig:euler_rmse_all}). While the initial assimilation reduces the error substantially for $\rho$, $\rho u$ and $\rho w$ of the EnDA run, the error increases over time and surpasses the error of the control EnNoDA run at approximately $t=225$ s. This increase in the error is due to the imbalances introduced by the chosen $(11 \times 11)$ grid point size of the localisation regions (more details are provided in section \ref{subsubsec:localisation_error} and Appendix B). For the EnDAB run, the imbalances are suppressed and the RMSEs are lower than those of the control EnNoDA run for all quantities over the entire simulation period.}

\subsubsection{Rising bubble}
\begin{figure*}[!ht]
    \centering
    \makebox[\textwidth][c]{
    \includegraphics[width=1.25\textwidth]{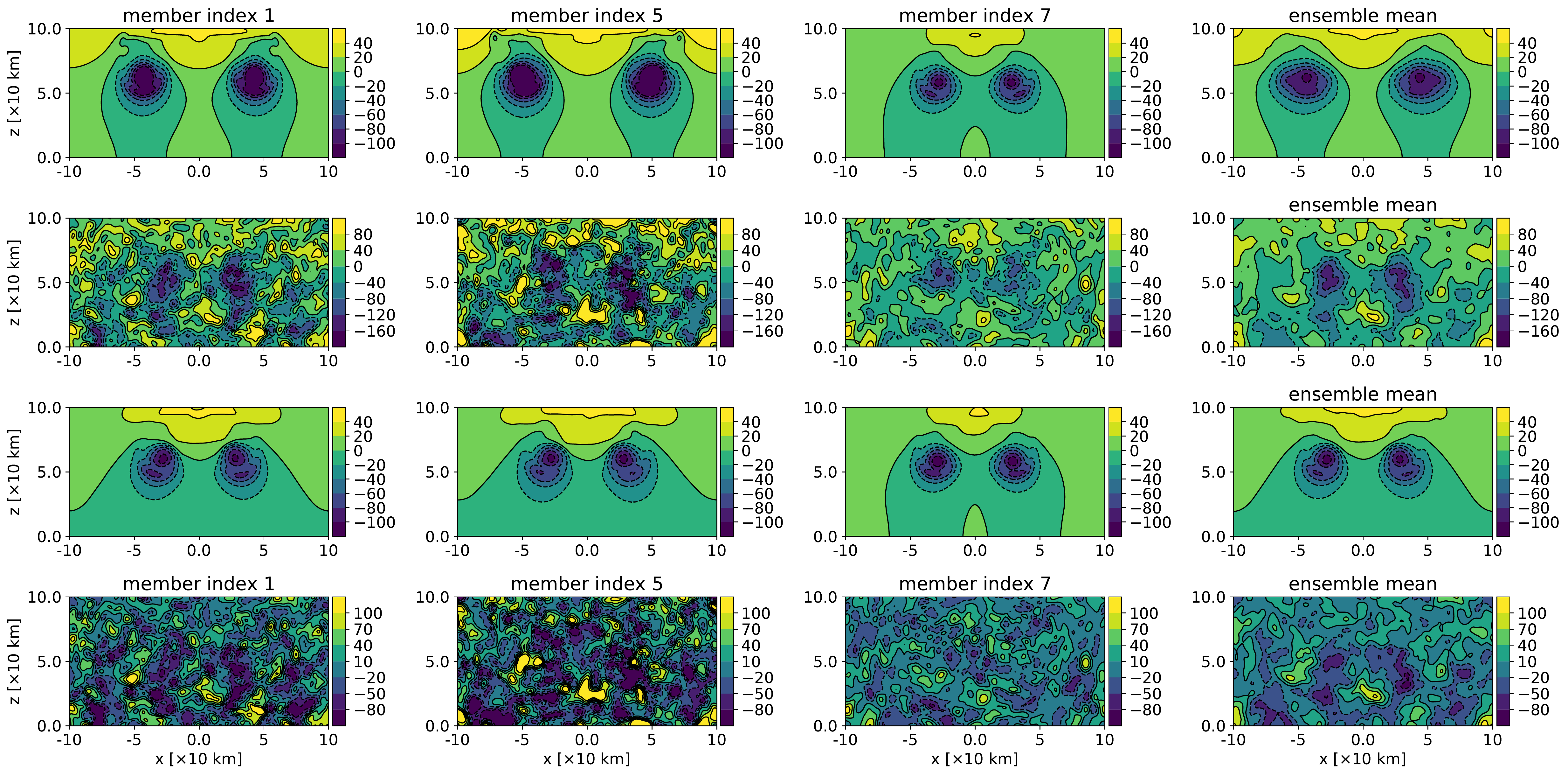}
    }
    \caption{Rising bubble: snapshots of pressure perturbation $p^\prime$. Ensemble members with index 1 (first column), 5 (second column) and 7 (third column) at $t=1000.0$ s with the ensemble mean (fourth column). Row-wise: EnNoDA run; contours in range $[-100,40]$ Pa with a $20$ Pa interval (first row), EnDA run; contours in range $[-160,80]$ Pa with a $40$ Pa interval (second row), EnDAB run; contours in range $[-100,40]$ Pa with a $20$ Pa interval (third row), and the difference between EnDA and EnDAB; contours in range $[-80,100]$ Pa with a $30$ Pa interval (fourth row). Negative contours are dashed.}
    \label{fig:rb_enses}
\end{figure*}

Figure \ref{fig:rb_enses} displays snapshots of pressure perturbation for the bubble case. In the EnNoDA run (first row) the bubbles in the ensemble attain different heights at the end of the simulation time and the ensemble mean is diffused, in line with the spread in the initial conditions used in generating the ensemble. Ensemble members with larger potential temperature perturbation rise faster. In the EnDA ensemble (second row), large-amplitude fast-mode imbalances are present while the ensemble mean of the bubble rotor \revision{positions} at the end time \revision{better approximates the true positions of the rotors}. For EnDAB (third row), the individual ensemble members are close to one another, as reflected in the ensemble mean. The ensemble \revision{better approximates} the truth and the fast-mode imbalances are suppressed. Moreover, the \revision{pressure footprints of the} bubble rotors are not visible \revision{in plots of the pressure differences} between the EnDA and EnDAB ensembles (fourth row), showing that the difference is \klein{predominantly} due to the presence of the imbalances only, and suggesting \klein{(last column)} that data assimilation is comparably effective in nudging the bubble towards the truth in both cases. \revision{Blending suppresses the imbalances while leaving the dynamics of the rising bubble largely unaffected}.

\begin{figure*}[ht]
    \centering
    \makebox[\textwidth][c]{
    \includegraphics[width=1.25\textwidth]{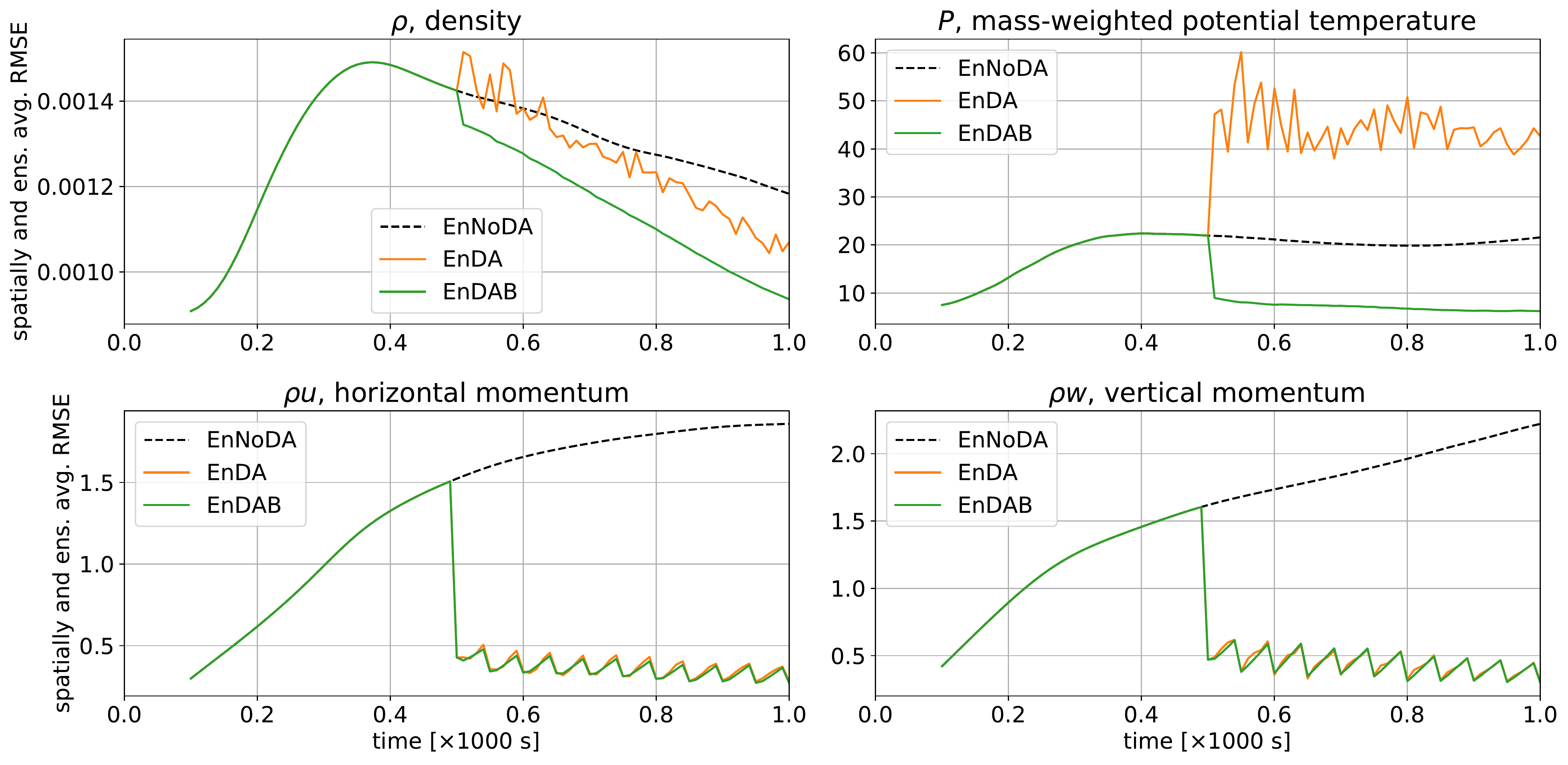}
    }
    \caption{Rising bubble: EnNoDA run (black), EnDA run (orange), EnDAB run (green). Assimilated quantities are $\rho u$ and $\rho w$. Spatially and ensemble averaged RMSE from $t=100.0$ s to $1000.0$ s for density $\rho$ (top left, [kg m$^{-2}$]), mass-weighted potential temperature $P$ (top right, [Pa]), and momenta $\rho u$, $\rho w$ (bottom left and right, [kg m$^{-1}$ s$^{-1}$]).}
    \label{fig:rb_rmse}
\end{figure*}

RMSE plots of data assimilation \revision{of the momentum fields} in the rising bubble experiment are \klein{shown} in Figure \ref{fig:rb_rmse}. The momentum fields are assimilated every 50.0 s after 500.0 s. This is visible in the momenta RMSE plots, where each downward step corresponds to \revision{one application of the} assimilation procedure. For EnDA, an error is introduced in the density $\rho$ and mass-weighted potential temperature $P$. Blending negates this and the EnDAB curves show a smooth profile, with RMSE lower than the control EnNoDA. As in the travelling vortex case, a jump is visible in the RMSE of $P$ at the first assimilation time for EnDA, and this corresponds to the imbalances introduced. See Appendix B on the scale analysis for more details. \revision{The ensemble spread and RMSE are again comparable in these runs (not shown).}

\subsubsection{Localisation region and imbalances}
\label{subsubsec:localisation_error}
\revision{In this section, results of the EnDA and EnDAB ensembles are investigated for varying localisation radii. Here the aim is not to obtain the optimal choice of the localisation radius but to illustrate its effect on the imbalances. All the quantities are assimilated for the travelling vortex test case and localisation regions of $(5 \times 5)$, $(21 \times 21)$ and $(41 \times 41)$ grid points are used. Otherwise, the setup follows the parameters laid out in sections \ref{subsubsec:vortex} and Table \ref{tab:experimental_setup}.

Increasing the size of the localisation region reduces the error for EnDA and EnDAB experiments. The error of the $(41 \times 41)$ travelling vortex EnDA experiment (cyan solid line in Figure \ref{fig:euler_locerrs}) is consistently lower than that of the control EnNoDA ensemble for all variables except for the pressure-related mass-weighted potential temperature $P$, for which the error jump is nevertheless small. This result suggests that the error in Figure \ref{fig:euler_rmse_all} for the EnDA ensemble with all quantities assimilated (orange dashed line) arises from the imbalances introduced through localisation. The amount of imbalances introduced can be reduced with large enough localisation regions. Imbalances are introduced by localisation even when taking into account a large proportion of the total grid points, e.g., $(41 \times 41)$ regions in a $(64 \times 64)$ mesh. Similar results are obtained for the rising bubble test (not shown).

\begin{figure*}[ht]
    \centering
    \makebox[\textwidth][c]{
    \includegraphics[width=1.25\textwidth]{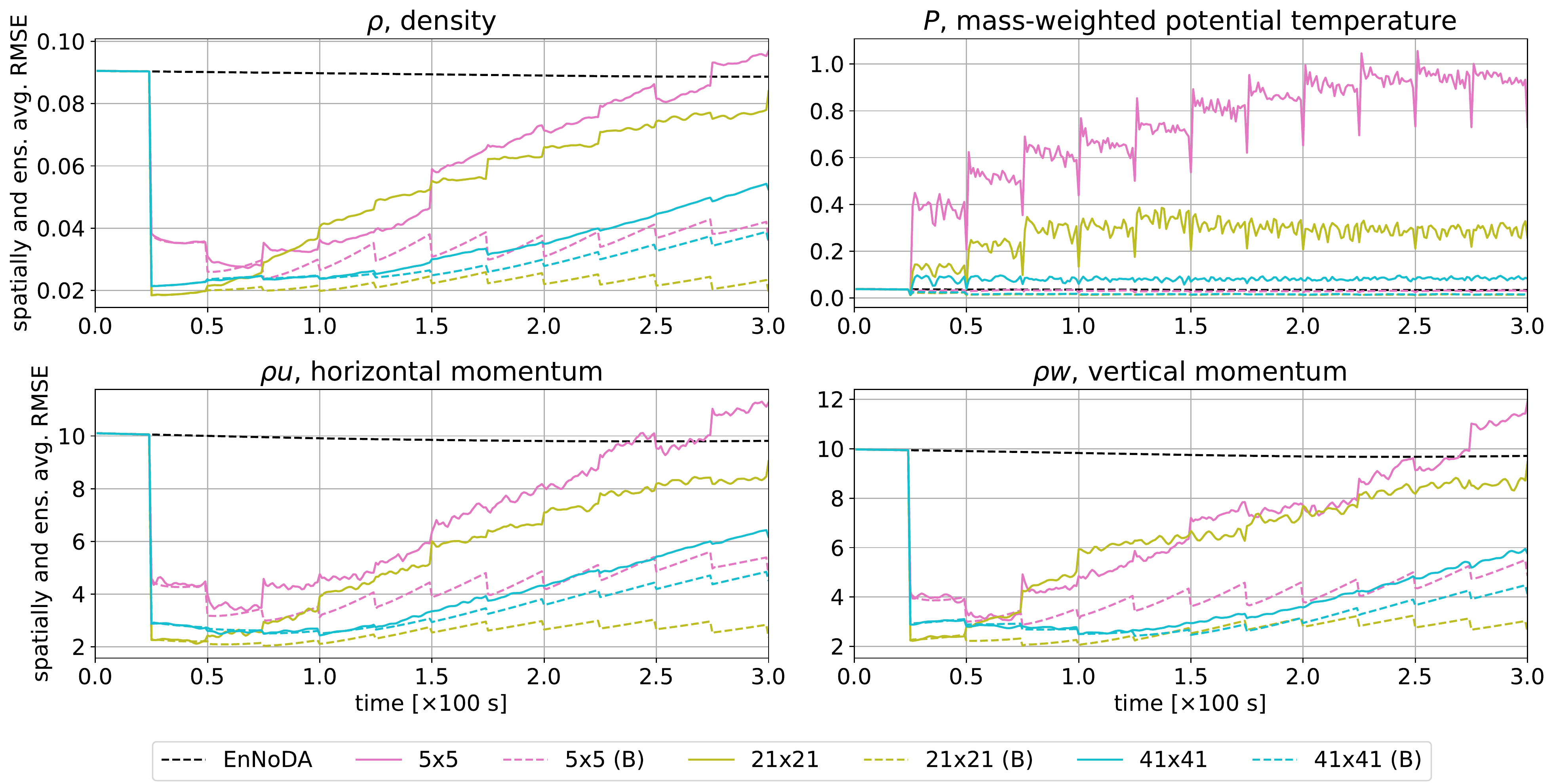}
    }
    \caption{\revision{Travelling vortex: EnNoDA (black) and experiments corresponding to localisation regions with $(5 \times 5)$ (magenta), $(21 \times 21)$ (yellow), and $(41 \times 41)$ (cyan) grid points. EnDA ensemble solutions are denoted with solid lines and EnDAB with dashed lines and with a `(B)' in the legend. Assimilated quantities are $\rho, \rho u, \rho w, P, \pi$. Spatially and ensemble averaged RMSE from $t=0.0$ s to $300.0$ s for density $\rho$ (top left, [kg m$^{-2}$]), mass-weighted potential temperature $P$ (top right, [kPa]), and momenta $\rho u$, $\rho w$ (bottom left and right, [kg m$^{-1}$ s$^{-1}$]). The RMSE of the initial ensemble is omitted.}}
    \label{fig:euler_locerrs}
\end{figure*}

For the localisation radii investigated in Figure \ref{fig:euler_locerrs}, the best performing EnDA solution performs worse than the worst performing EnDAB solution (compare the solid cyan and dashed magenta lines). In the case of the EnDA solutions (solid lines), the fast-mode imbalances introduced are a significant source of error, see section \ref{subsubsec:vortex} for more details. As a consequence, any reduction in the imbalances introduced will lead to a corresponding decrease in the error. This corroborates the decrease in the RMSE for increasing localisation sizes. For the EnDAB solutions (dashed lines), the balanced structure is preserved and there is an optimal localisation length scale. Too small a localisation region will lead to an under-sampling of the vortex dynamics, while too large a localisation region may distort the vortex dynamics by oversampling the background dynamics. This is observed in the larger RMSE for the $(41 \times 41)$ solution (dashed cyan) when compared to the $(21 \times 21)$ solution (dashed yellow). These results indicate that, for blended data assimilation of the travelling vortex experiment, a moderate localisation region of approximately a third of the grid-size, i.e. $(21 \times 21)$ cells in a $(64 \times 64)$ grid, yields the optimal analysis fields.
}

\section{Discussion and conclusion}
\label{sec:diss_conc}
This paper has presented a new conceptual framework for balanced data assimilation based on blended numerical models. Using a \klein{discrete time-level numerical analysis for the Exner pressure field} and a careful choice of pressure perturbation variables, the blended soundproof-compressible modelling framework of \cite{bok2014} has been substantially upgraded by a functionality to switch between equation sets in a single time-step.

\revision{In idealised numerical experiments with a travelling vortex and a gravity-driven warm air bubble, a single time-step in the pseudo-incompressible limit regime was sufficient to recover a balanced state starting from imbalanced initial data. Moreover, the blended model yielded leftover acoustics with amplitude more than one order of magnitude smaller than the ones generated at the onset with the fully compressible model.} The amplitude reduction is a sizeable improvement over the scores of \cite{bok2014} who, in addition, \klein{needed several time steps in a hybrid soundproof-compressible configuration with non-integer values of the blending parameter $\alpha_P$ to achieve their best level of noise reduction}.

The upgraded blended model has then been combined with a data assimilation engine and deployed as a tool to reduce imbalances introduced by regular assimilation of data within model runs. Numerical results on ensemble data assimilation with \revision{and without} blending showed that while data assimilation alone produced imbalances that effectively destroyed \revision{important qualitative features of the solution in one of the test cases, data assimilation together with blending strongly reduced those imbalances and lead to recovery of} accurate results. \revision{Moreover, blended data assimilation was effective despite the untuned data assimilation parameters used in the investigations. Throughout our study,} a single time-step spent in the pseudo-incompressible limit regime after the assimilation of data was sufficient to restore the balanced state, as \revision{documented by strongly} reduced RMSEs with the blended model.

For ensemble data assimilation experiments with the travelling vortex, assimilation of the momentum fields alone was found to be insufficient. Over longer simulations, the ensemble with balanced data assimilation carried larger errors than the control ensemble without data assimilation. See the green solid curves in Figure \ref{fig:euler_rmse_all}. \klein{We have traced the origin of this result back to an issue of controllability \citep{Jazwinski2007}: This test case involves large, dynamically relevant potential temperature variations whose deviation from the truth cannot be \revision{controlled} at all when only the momentum field is assimilated. In fact, a test with an analogous vortex that has constant entropy initial data yields results \revision{(not shown)} close in quality to those of the rising thermal test when only momentum is assimilated. The issue could then be solved by assimilation of all variables. F}urther investigation is warranted on how the effectiveness of data assimilation can be improved under such circumstances \klein{without the need to observe all state variables}. A scale analysis (Appendix B) corroborated the insight that the RMSE increase introduced by the assimilation of data corresponds to the fast-mode imbalances seen in the plots of the individual ensemble members and the ensemble mean. \revision{In this sense, our experiments make a case for investigations involving relatively simple idealised test cases, as we were able to gain some analytical understanding of the sources and consequences of errors and imbalances. \revision{Nevertheless, further studies based on more realistic scenarios will be required to demonstrate that the presented approach and its extensions will actually enable quantifiable improvements of numerical weather prediction skill scores}}.

\revision{In the experiments involving ensemble data assimilation with different localisation radii, blended data assimilation yielded, for all localisation sizes, substantial improvements to the RMSE relative to the plain data assimilation without a balancing procedure. In fact, the best-performing data assimilation-only run still produced worse results than the worst-performing run with blending. Furthermore, the recovery of a balanced vortex structure by blended data assimilation turned out to be sensitive to the choice of localisation radius, with best results obtained at some intermediate size of the localization domains. In contrast, increasing the localisation size for an ensemble with data assimilation without blending was sufficient to decrease its error. Yet, since the imbalances by far dominate the overall error in this case, this is just a reflection of the expected noise reduction resulting from increased smoothing of the assimilated information.}

\revision{In numerical weather prediction, methods to damp or remove acoustic imbalances have long been employed \citep[e.g. ][]{daley1988normal, skamarock1992stability, dudhia1995reply, klemp2018damping}.  Moreover, practical application of sequential data assimilation procedures will generally excite all rapidly oscillatory modes of the compressible system, and filtering techniques are used to negate these unphysical imbalances \citep{ha2017ensemble}. In this context, the results presented in this paper are encouraging in that blended data assimilation was able to suppress acoustic noise and recover balanced analysis fields, albeit for idealised test cases. To the best of the authors' knowledge, this is the first study of a dynamics-driven method to suppress acoustic noise arising from the sequential assimilation of data.}

In addition, the results presented in this paper \klein{prepare the ground} for future work in a number of areas. \revision{In general, the performance of a data assimilation method can be improved by tuning its adjustable parameters. Here, however, we consciously employed an untuned data assimilation scheme known to produced unphysical imbalances to test the efficacy of our dynamics-driven method in removing them. Consequently, a comprehensive study similar to \cite{popov2019bayesian} on multivariate tuning of the LETKF and localisation parameters for the blended numerical model will be an avenue for future improvements of our approach. The study could also compare our method with existing balancing strategies, e.g., the IAU and the DFI, following \cite{polavarapu2004relationship}. To ensure a fair comparison, optimisations of the IAU along the lines of \cite{lei2016four} and \cite{he2020impacts} may have to be carried out. A comparison of the effects of our dynamics-driven method on the slower dynamics against those of the DFI and IAU, which act as low-pass filters \citep{houtekamer2016review, polavarapu2004relationship}, will be particularly insightful.}

\revision{Despite the untuned data assimilation scheme used, the blended model has given promising results, although thus far only for idealised test cases\label{idealised_3}. Another  natural evolution will hence involve model performance on more realistic three-dimensional moist dynamics scenarios with bottom topography \citep{ONeillKlein2014,duarte2015low} and on benchmarks at larger scales \citep{skamarock1994efficiency, bk2019}}.

\revision{Although presented and refined here for the blending between the compressible Euler equations and the pseudo-incompressible model only, the methodology translates to other scenarios as long as one can formulate the according projection onto appropriate reduced dynamics via the elliptic pressure correction. Models imposing a divergence constraint on the weighted velocity field as well as frameworks blending between nonhydrostatic and hydrostatic dynamics will naturally fit into the present approach.
}

\revision{Specifically, t}he numerical scheme proposed by \cite{bk2019} enables solution of the hydrostatic system in the large-scale limit in addition to the small-scale low Mach number limit considered in this paper. Therefore, a blended data assimilation framework such as the one presented here could be enhanced with hydrostatic blending and used in a two-way blended pseudo-incompressible / hydrostatic / compressible model \citep{klein2016} \revision{exploiting the different dynamics in the equation sets}. 

\revision{Moreover, the theoretical framework developed in that paper also included the unified model by \citet{arakawa2009unification} as one of the reduced models. Thus, after an appropriate extension of the present numerical scheme, yet another framework for blended data assimilation can be developed. In fact, a variant of the fully compressible/Arakawa-Konor model pair has recently been presented by \cite{qaddouri2021implementation}, and a related blending approach will allow for the filtering of smaller-scale acoustic noise while leaving the Lamb wave components dynamically unaffected.}
Investigations similar to the ones in this paper can then be made on balancing initial states and data assimilation for small- to planetary-scale dynamics using the resulting doubly blended model framework. Internal waves play an important role for atmospheric dynamics and they should not be removed indiscriminately after a data assimilation step. Therefore, the identification and removal of unwanted internal wave noise while keeping the physically meaningful wave spectrum is an additional challenge that will require further theoretical developments beyond the scope of this paper.

\revision{More generally, semi-implicit compressible models feature in several dynamical cores used by weather centres worldwide. Notable examples include the currently operational hydrostatic IFS spectral transform model in use at the European Centre for Medium Range Weather Forecasts \citep[ECMWF,][]{wedi2013fast}, and the Met Office's Unified Model \citep{davies2005new, wood2014inherently}, which has a hydrostatic-nonhydrostatic switch. ECMWF's next-generation nonhydrostatic compressible dynamical core, IFS-FVM \citep{kuhnlein2019fvm}, actually uses a numerical discretisation akin to the one considered in this paper and would therefore be an ideal candidate for a first implementation of the blended tools in a semi-operational model.  In addition, our approach will bear particular relevance to fully compressible operational models featuring the option of selectively employing the dynamics of a limit model \citep{wood2014inherently, melvin2019mixed, voitus2019, qaddouri2021implementation}.}

In this context, multimodel numerics with seamless switching could contribute to creating a level playing field to evaluate accuracy and performance with different equation sets in \revision{the same} dynamical core. The positive evidence provided here in balancing data assimilation shows\klein{, in the authors' view, a considerable potential and potential impact of deploying the blended model} framework across the whole forecast model chain.

\paragraph{Acknowledgments}
R.C., G.H. and R.K. thank the Deutsche Forschungsgemeinschaft for the funding through the Collaborative Research Center (CRC) 1114 ``Scaling cascades in complex systems'', \klein{Project Number 235221301,}  Project A02: ``Multiscale data and asymptotic model assimilation for atmospheric flows''. T.B. was supported by the ESCAPE-2 project, European Union’s Horizon 2020 research and innovation programme (grant agreement No 800897).


\renewcommand{\thefigure}{A\arabic{figure}}
\setcounter{figure}{0}

\begin{appendices}
\section{LETKF Algorithm}
\label{apx:LETKF}
The Local Ensemble Transform Kalman Filter (LETKF) algorithm presented here is a summary of the algorithm published by \cite{hunt2007} in their paper, adapted to the blended numerical framework.

Start with an ensemble of $K$ state vectors, $\{ \mathbf{x}_{k,[g]}^{f} \} \in \mathbb{R}^{m_{[g]}}$ for $k=1,\dots,K$. Furthermore, assume that a set of observations $\mathbf{y}_{\text{obs}, [g]} \in \mathbb{R}^{l_{[g]}}$ with a known covariance $\pmb{\mathsf{R}}_{[g]} \in \mathbb{R}^{l_{[g]} \times l_{[g]}}$ is available. Here, $m$ and $l$ represent the dimension of the state and observation spaces and the subscript $\mathit{[g]}$ represents the global state space, \emph{i.e.} localisation has not been applied.

\begin{enumerate}
\item Apply the forward operator $\mathcal{H}$ to obtain the state vectors in the observation space,
\begin{equation}
\mathcal{H}(\mathbf{x}_{k,[g]}^{f}) = \mathbf{y}_{k,[g]}^{f} \in \mathbb{R}^{l_{[g]}}.
\end{equation}
\item Stack the anomaly of the state and observation vectors to form the matrices,
\begin{align}
\pmb{\mathsf{X}}^f_{[g]} &= \left[\mathbf{x}^{f}_{1,[g]} - \bar{\mathbf{x}}_{[g]} \left\vert ~ \dots ~ \right\vert \mathbf{x}^{f}_{K,[g]} - \bar{\mathbf{x}}_{[g]} \right] ~ \in \mathbb{R}^{m_{[g]} \times K},\label{eqn:local_X}
\\
\pmb{\mathbf{Y}}^f_{[g]} &= \left[\mathbf{y}^{f}_{1,[g]} - \bar{\mathbf{y}}_{[g]} \left\vert ~ \dots ~ \right\vert \mathbf{y}^{f}_{K,[g]} - \bar{\mathbf{y}}_{[g]}\right] ~ \in \mathbb{R}^{l_{[g]} \times K},\label{eqn:local_Y}
\end{align}
where $\bar{\mathbf{x}}_{[g]}$ ($\bar{\mathbf{y}}_{[g]}$) is the mean of the state vectors (in observation space) over the ensemble, e.g.
\begin{equation}
\bar{\mathbf{x}}_{[g]} = \frac{1}{K} \sum_{k=1}^K \mathbf{x}^{f}_{k,[g]} ~ \in \mathbb{R}^{m_{[g]}}.
\end{equation}
\item From $\pmb{\mathsf{X}}^f_{[g]}$ and $\pmb{\mathsf{Y}}^f_{[g]}$, select the local $\pmb{\mathsf{X}}^f$ and $\pmb{\mathsf{Y}}^f$.

\item From the global observations $\mathbf{y}_{\text{obs}, [g]}$ and observation covariance $\pmb{\mathsf{R}}_{[g]}$, select the corresponding local counterparts $\mathbf{y}_{\text{obs}}$ and $\pmb{\mathsf{R}}$. Notice that the subscript $\mathit{[g]}$ is dropped when representing the local counterparts.

\item Solve the linear system $\pmb{\mathsf{R}} \pmb{\mathsf{C}}^T = \pmb{\mathsf{Y}}^f$ for $\pmb{\mathsf{C}} \in \mathbb{R}^{K \times l}$.

\item Optionally, apply a localisation function to $\pmb{\mathsf{C}}$ to modify the influence of the surrounding observations.

\item Compute the $K \times K$ gain matrix,
\begin{equation}
{\pmb{\mathsf{K}}} = \left[ (K-1) \frac{\pmb{\mathsf{I}}}{b} + \pmb{\mathsf{C}} \pmb{\mathsf{Y}}^f \right]^{-1},
\label{eqn:ens_covar}
\end{equation}
where $b > 1$ is the ensemble inflation factor.

\item Compute the $K \times K$ analysis weight matrix,
\begin{equation}
\pmb{\mathsf{W}}^a = \left[ (K-1) ~ {\pmb{\mathsf{K}}} \right]^{1/2}.
\end{equation} 
\item Compute the $K$-dimension vector encoding the distance of the observations from the forecast ensemble
\begin{equation}
\bar{\mathbf{w}}^a = {\pmb{\mathsf{K}}} \pmb{\mathsf{C}} \left( \mathbf{y}_{\rm obs} - \bar{\mathbf{y}}^f \right),
\end{equation}
and add $\bar{\mathbf{w}}^a$ to each column of $\pmb{\mathsf{W}}^a$ to get a set of $K$ weight vectors $\{ \mathbf{w}^{a}_k \}$ with $k=1,\dots,K$.
\item From the set of weight vectors, compute the analysis for each ensemble member,
\begin{equation}
\mathbf{x}^{a}_k = \pmb{\mathsf{X}}^f \mathbf{w}^{a}_k + \bar{\mathbf{x}}^f, \quad \text{for } k=1,\dots,K.
\end{equation}
\item Finally, recover the global analysis ensemble $\{ \mathbf{x}_{k,[g]}^{a} \}$, $k=1,\dots,K$.

This recovery depends on how the local regions were selected in \eqref{eqn:local_X} and \eqref{eqn:local_Y}. For local region surrounding the grid point under analysis, the reassembly of the global analysis ensemble is done by reassembling the analysed grid points back into the global grid.
\end{enumerate}

\renewcommand{\thefigure}{B\arabic{figure}}
\setcounter{figure}{0}
\section{Scale analysis for the\\data assimilation error in the pressure-related fields}
\label{apx:scale_analysis}
Figures \ref{fig:euler_rmse_all} and \ref{fig:rb_rmse} show that the assimilation of only the momentum fields leads to a jump in RMSE in the non-momentum fields\revision{, and the assimilation of all quantities in Figure \ref{fig:euler_rmse_all} leads to a jump in RMSE in the pressure-related $P$ field}. This increase in the error occurs after the first assimilation time and remains of the same order of magnitude for the duration of the simulation, quantifying the imbalance introduced by data assimilation. The imbalance can be characterised by a scale analysis \citep{klein2001asymptotic}.

The assimilation of the momentum fields leads to a change in the divergence of the velocity fields,
\begin{equation}
    \nabla \cdot (\delta \mathbf{v}) = \frac{\partial \delta u}{\partial x} + \frac{\partial \delta w}{\partial z},
    \label{eqn:change_in_div_velocity}
\end{equation}
where $(\delta u, \delta w)$ are the changes in the velocity fields due to the assimilation of momenta in the vertical slice experiments. \eqref{eqn:change_in_div_velocity} has the units [s$^{-1}$].

Observe from Figures \ref{fig:euler_enses_3.0} and \ref{fig:rb_enses} that the imbalance introduced by data assimilation are fast-mode acoustic waves. This effect is modelled as a wave oscillating with the peak amplitude right after the assimilation of data at the grid point under analysis. Therefore, for an oscillating wave excited at grid point $(x_i,z_j)$, the maximum amplitude of the imbalances is
\begin{align}
    (\nabla \cdot \delta \mathbf{v})_{(i,j)} \int_{0}^{t_{\text{ac}}} \cos \left(\frac{\pi}{2} \frac{t}{t_{\text{ac}}} \right) dt &= \frac{2 t_{\text{ac}}}{\pi}(\nabla \cdot \delta \mathbf{v})_{(i,j)} \left[ \int_{0}^{\pi/2} \cos({\xi}) d \xi  \right]  \nonumber \\
    & = \frac{2 t_{\text{ac}}}{\pi}(\nabla \cdot \delta \mathbf{v})_{(i,j)},
    \label{eqn:oscillating_imbalance}
\end{align}%
The acoustic timescale $t_{\text{ac}}$ is chosen as the timescale of the largest perturbations introduced. This is the time a wave takes to traverse to the edge of the $(11 \times 11)$ grid points local region from the analysis grid point. Therefore,
\begin{equation}
    t_{\text{ac}} = \frac{11}{2} \frac{dx}{c_{\text{ref}}},
    \label{eqn:t_ac}
\end{equation}
where $dx$ is the constant grid-size and $c_{\text{ref}}$ the speed of sound. \eqref{eqn:t_ac} has units [s] and \eqref{eqn:oscillating_imbalance} is dimensionless.

As $p = \rho c_{\text{ref}}^2$, the contribution to the pressure from $\nabla \cdot (\delta \mathbf{v})_{(i,j)}$ is computed by
\begin{equation}
    \frac{2 t_{\text{ac}}}{\pi} \nabla \cdot (\delta \mathbf{v})_{(i,j)} \rho_{(i,j)} c_{\text{ref}}^2 \sim \hat{p}_{(i,j)},
\end{equation}
which has the units of [Pa]. The hat \,$\mathit{\hat{~}}$\, signifies that the quantity is obtained from scale analysis. Finally, use the equation of state \eqref{eqn:eos} to obtain an estimate for $\hat{P}$.

\begin{figure*}[!ht]
    \centering
    \makebox[1.0\textwidth][c]{
    \begin{tabularx}{1.25\textwidth}{cc}
    \includegraphics[width=0.6\columnwidth]{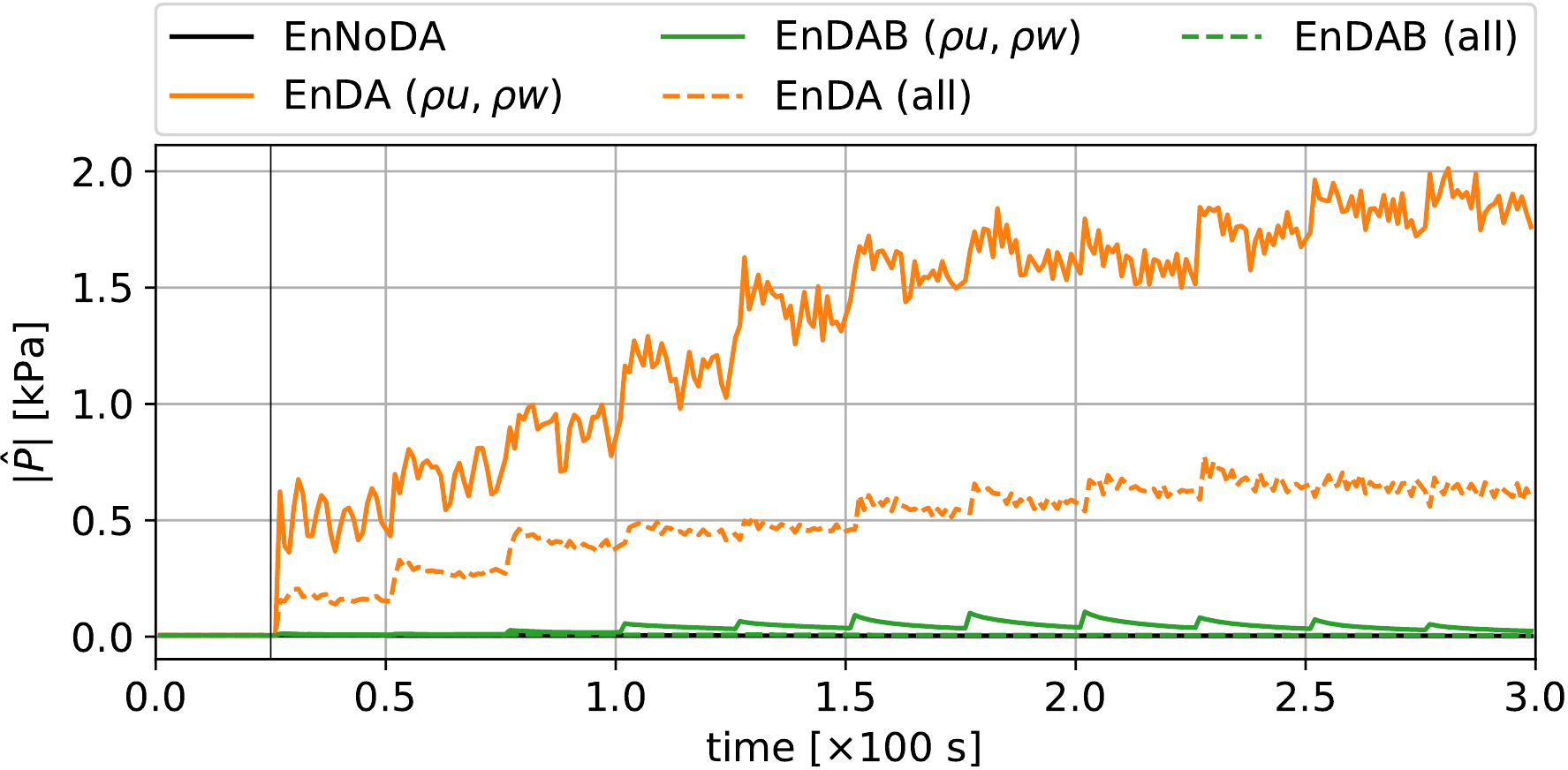}
    &
    \includegraphics[width=0.6\columnwidth]{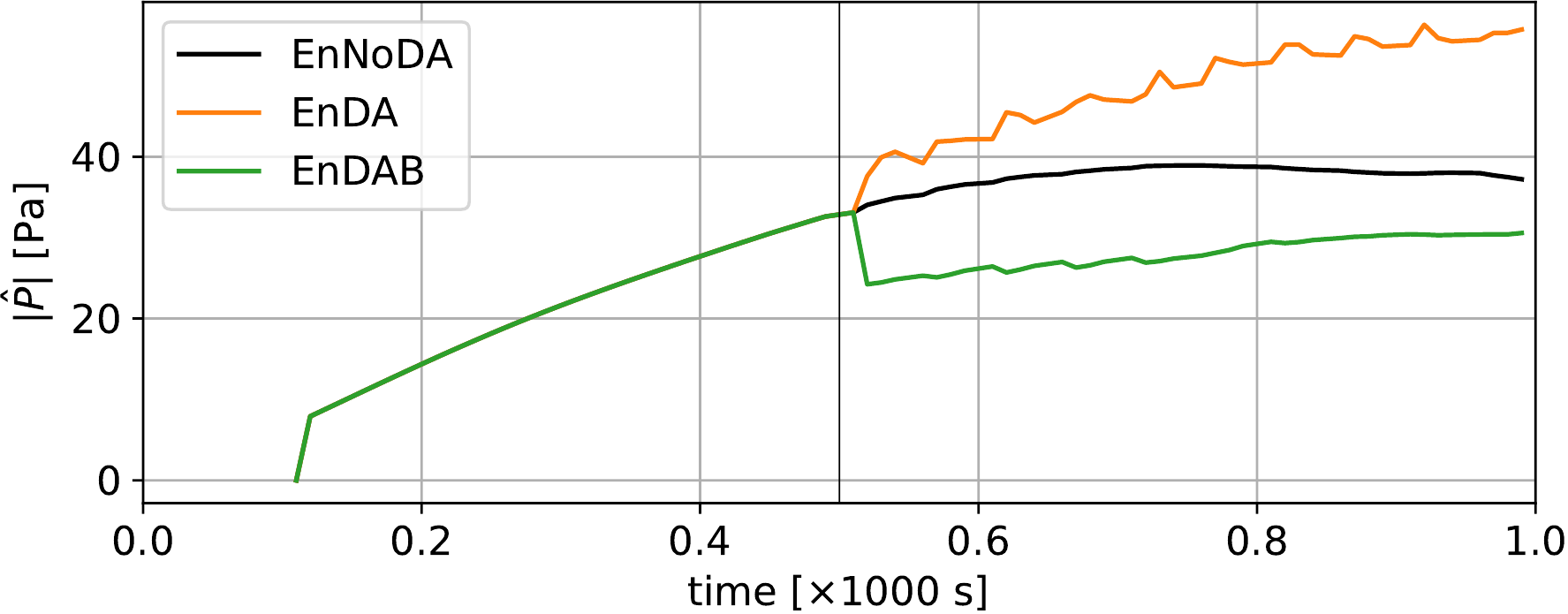}
    \end{tabularx}
    }
    \caption{Scale analysis of the contribution to the mass-weighted potential temperature $P$ from the divergence of the velocity fields for the travelling vortex ensemble (left) and the rising bubble ensemble (right). \revision{In the legend, (all) represents the travelling vortex ensembles with all quantities, $\rho, \rho u, \rho w, P, \pi$, assimilated.} The first assimilation time is marked with a vertical solid black line.}
    \label{fig:scale_analysis}
\end{figure*}

For comparison with the RMSE, the norm is taken for $\hat{P}$, given by
\begin{equation}
    \left| \, \hat{P} \, \right| = \frac{1}{K} \sum_{k}^{K} \left[ \sqrt{\frac{1}{N_x \times N_z} \sum_{i,j}^{N_x,N_z} \left( \hat{P}_{(i,j)} \right)^2} ~ \right]_k,
\end{equation}
where $k$ indexes the $K$ ensemble members and $N_x$ and $N_z$ are the number of grid points in the $x$ and $z$ coordinates.

Figure B1 shows the results of scale analysis for the two test cases. Results at assimilation time are omitted. Scale analysis yields EnDA results for \revision{$\left| \, \hat{P} \, \right|$} that are of the same order of magnitude as the jumps in the RMSE plots (Figures \ref{fig:euler_rmse_all} and \ref{fig:rb_rmse}) with a similar profile over time. Scale analysis characterises the error jump in the thermodynamical RMSE plots as fast-mode imbalances introduced through data assimilation.

\end{appendices}

\bibliographystyle{abbrvnat}
\bibliography{references}

\end{document}